\documentclass[12pt,twoside]{amsart}
\usepackage{amssymb,amsmath,amsthm, amscd, enumerate, mathrsfs}
\usepackage{graphicx, hhline}
\usepackage[all]{xy}
\usepackage[shortlabels]{enumitem}
\usepackage{mathtools}
\mathtoolsset{showonlyrefs}
\usepackage{color}
\usepackage{hyperref}

\title[Remarks on minimal model theory]{Remarks on the minimal model theory for log surfaces in the analytic setting}
\author{Nao Moriyama} 
\address{Department of 
Mathematics, Graduate School of Science, 
Kyoto University, Kyoto 606-8502, Japan}
\email{moriyama.nao.22s@st.kyoto-u.ac.jp}
\keywords{Minimal model program}
\subjclass[2020]{Primary 14E30; Secondary 32C15}

\newtheorem{theorem}{Theorem}[section]
\newtheorem{lemma}[theorem]{Lemma}
\newtheorem{corollary}[theorem]{Corollary}

\theoremstyle{definition}
\newtheorem{definition}[theorem]{Definition}
\newtheorem{remark}[theorem]{Remark}
\newtheorem*{ack}{Acknowledgments}  
\newtheorem{step}{Step}
\newtheorem{case}{Case}

\makeatletter

\@addtoreset{equation}{section}
\makeatother
\newcounter{stepcounter}
\renewcommand{\step}{%
  \stepcounter{stepcounter}% 
  \medskip\noindent\textbf{Step \arabic{stepcounter}. }
}

\begin{document}

\begin{abstract}
We discuss the relative log minimal model theory for log surfaces in the analytic setting. More precisely, we show that the minimal model program, the abundance theorem, and the finite generation of log canonical rings hold for log pairs of complex surfaces which are projective over complex analytic varieties.
\end{abstract}

\maketitle

\tableofcontents

\section{Introduction}
We study the minimal model theory for log surfaces in the analytic case. In order to extend the theory from the algebraic setting to the analytic setting, we introduce several specific techniques.

Before stating the main result, we describe the conditions adopted in the analytic case.
We say that $(X, Y, W, \pi, \Delta)$ satisfies the condition ($\bigstar$) if it satisfies the following conditions.
\begin{itemize}
\item $(X, \Delta)$ is a pair of a normal complex surface $X$ and a boundary $\mathbb{R}$-divisor $\Delta$ on $X$. 
\item $K_X+\Delta$ is $\mathbb{R}$-Cartier.
\item $\pi\colon X\rightarrow Y$ is a projective surjective morphism of complex analytic varieties.
\item $W$ is a compact subset of $Y$ such that $\rho(X/Y;W)<\infty$.
\end{itemize} 

In the analytic case, the notation $K_X$ is somewhat delicate. This is because a divisor corresponding to the canonical line bundle may not be globally well-defined. Therefore, the second condition should more precisely be stated as: $K_X+\Delta$ is $\mathbb{R}$-Cartier at every point $x\in X$. In this paper, we adopt the simplified notation without further clarification for the sake of simplicity.
The last condition is also specific to the analytic setting. We replace $Y$ with a small open neighborhood of $W$ if necessary.

For later use, we collect here four cases.
\begin{enumerate}[{Case} (\Alph*)]
\item $K_X$ is $\mathbb{Q}$-Cartier and every prime divisor on $X$ which is mapped into $W$ is $\mathbb{Q}$-Cartier.
\item $X$ is $\mathbb{Q}$-factorial over $W$.
\item $(X, \Delta)$ is log canonical.
\item $X$ has only rational singularities.
\end{enumerate}
In what follows, we study the minimal model theory in each of Cases (A)–(D). Case (A) is a new setting.
The following theorem is the main result of this paper.
\begin{theorem}[{see, \cite[Theorem 1.1]{Fjn12}}] \label{Theorem 1.1}
Assume that $(X, Y, W, \pi, \Delta)$ satisfies the condition \textnormal{($\bigstar$)}
and that one of \textnormal{Cases (A)–(D)} holds.
We shrink $Y$ around $W$ suitably.
Then there is a sequence of at most $\rho(X/Y;W)-1$ contractions
$$
(X, \Delta)=(X_0, \Delta_0)\stackrel{\varphi_0}{\rightarrow}(X_1, \Delta_1)\stackrel{\varphi_1}{\rightarrow}\cdots\stackrel{\varphi_{k-1}}{\rightarrow}(X_k, \Delta_k)=(X^*, \Delta^*)
$$
over $Y$ such that
\begin{itemize}
\item \textnormal{(Good minimal model)} if $K_X+\Delta$ is $\pi$-pseudo-effective, then $K_{X^*}+\Delta^*$ is semi-ample over some open neighborhood of $W$, and
\item \textnormal{(Mori fiber space)} if $K_X+\Delta$ is not $\pi$-pseudo-effective, then $(X^*, \Delta^*)$ is a Mori fiber space.
\end{itemize}
We note that $K_{X_i}$ is $\mathbb{Q}$-Cartier and every prime divisor on $X_i$ which is mapped into $W$ is $\mathbb{Q}$-Cartier for every $i$ in \textnormal{Case (A)}, that $X_i$ is $\mathbb{Q}$-factorial over $W$ for every $i$ in \textnormal{Case (B)}, that $(X_i, \Delta_i)$ is log canonical for every $i$ in \textnormal{Case (C)},
and that $X_i$ has only rational singularities for every $i$ in \textnormal{Case (D)}.
\end{theorem}

\begin{remark}
When $\operatorname{dim}X>2$, each contraction morphism given by the minimal model program exists only after shrinking around $W$ as mentioned in \cite[Theorem 1.7]{Fjn22}. Then we have to repeatedly replace $Y$ with a small open neighborhood of $W$.
In contrast, when $\operatorname{dim}X=2$, we can construct an open neighborhood $U$ of $W$ on $Y$ such that, once we replace $Y$ with $U$, each step $\varphi_i$ in the above theorem exists.
We note that this $U$ depends only on the compact subset $W$ and the morphism $\pi\colon X\rightarrow Y$, and it is independent of the choice of the divisor $\Delta$.
\end{remark}

\begin{definition}[Mori fiber space]
Assume that $(X, Y, W, \pi, \Delta)$ satisfies the condition ($\bigstar$).
Let $\varphi\colon X\rightarrow Z$ be a projective morphism of normal complex analytic varieties over $Y$. Then $\varphi\colon(X, \Delta)\rightarrow Z$ is a \textit{Mori Fiber space} over $Y$ if
\begin{enumerate}
    \item[(i)] $\varphi$ is a contraction morphism associated to a $(K_X+\Delta)$-negative extremal ray of $\overline{\text{NE}}(X/Y; W)$, and
    \item[(ii)] $\operatorname{dim}Z<2$.
\end{enumerate}
\end{definition}

In the algebraic case, Fujino showed that the analogues of Theorem \ref{Theorem 1.1} in Cases (B) and (C) are true over $\mathbb{C}$ (cf.~\cite{Fjn12}).
In positive characteristics, Tanaka proved that the corresponding statements hold (cf.~\cite{Tnk14}).

We cannot directly use this algebraic framework, due to an issue in the analytic setting.
A key difficulty is that $\mathbb{Q}$-factoriality does not behave well.
We take a compact subset $W'$ of $Y$ with $W'\subset W$. Then, $X$ is not necessarily $\mathbb{Q}$-factorial over $W'$ even if $X$ is $\mathbb{Q}$-factorial over $W$. This is because there may exist a divisor defined over an open neighborhood of $\pi^{-1}W'$ which can not be extended to a divisor defined over an open neighborhood of $\pi^{-1}W$ (cf.~\cite[Remark 2.39]{Fjn22}). 
This fact is troublesome because we sometimes want to take a Stein open subset $U$ of $Y$, which may not contain $W$, and replace $Y$ with $U$ in the analytic case. 
In order to address this problem, we replace the assumption of $\mathbb{Q}$-factoriality with a weaker condition that is local on $Y$. Specifically, we assume that $K_X$ is $\mathbb{Q}$-Cartier and that every prime divisor on $X$ mapped into $W$ is $\mathbb{Q}$-Cartier; this is precisely the condition stated in \textnormal{Case (A)}.
Moreover, we can reduce Case (B), where $X$ is $\mathbb{Q}$-factorial over $W$, to Case (A) in order to prove the main result. Indeed, if $X$ is $\mathbb{Q}$-factorial over $W$, then every prime divisor on $X$ mapped into $W$ is $\mathbb{Q}$-Cartier. In addition, by suitably shrinking $Y$ around $W$, we may also assume that $K_X$ is $\mathbb{Q}$-Cartier.

Theorem \ref{Theorem 1.1} is a consequence of the following two theorems: the minimal model program and the abundance theorem.
In the proof of the minimal model program in Case (A), we use the condition that every prime divisor on $X$ mapped into $W$ is $\mathbb{Q}$-Cartier, but we need not require $K_X$ to be $\mathbb{Q}$-Cartier. The latter condition is used in the proof of the following abundance theorem in Case (A), which constitutes the most technically challenging part of this paper.
\begin{theorem}[{Abundance theorem, see Theorem \ref{thm: abundance for R-divisors}}]
Assume that $(X, Y, W, \pi, \Delta)$ satisfies the condition \textnormal{($\bigstar$)}
and that one of \textnormal{Cases (A)–(D)} holds.
If $K_X+\Delta$ is nef over $W$, then it is semi-ample over a neighborhood of $W$.
\end{theorem}

We encounter two additional difficulties in the analytic setting. First, we need to prove the openness of nefness. In the algebraic case, it suffices to show that if $K_X + \Delta$ is nef, then it is semi-ample. Thus, the openness of nefness is not a problem. In contrast, in the analytic setting, we have to show that if $K_X + \Delta$ is nef over $W$, then it is nef even over a neighborhood of $W$. In this paper, we prove this openness of nefness, in the case where $X$ is a complex surface (cf.~Lemma \ref{lem: openness of nef} below). Second, the technique of compactification is not available. In the algebraic setting, Fujino reduced the relative abundance theorem to the absolute case (i.e.,~$\dim Y = 0$) by using compactification (cf.~\cite{Fjn12}). However, in the analytic setting, such a reduction is not possible. Therefore, we need to prove the relative abundance theorem directly in the cases where $\dim Y = 1$ and $\dim Y = 2$, respectively.\\
\par
As a corollary of the abundance theorem, we obtain the following finite generation of log canonical rings.
\begin{corollary}[{see, \cite[Theorem 4.5]{Fjn12}}] \label{Corollary 1.5}
Assume that $(X, Y, W, \pi, \Delta)$ satisfies the condition \textnormal{($\bigstar$)}, $\Delta$ is a $\mathbb{Q}$-divisor and that one of \textnormal{Cases (A) and (B)} holds.
After shrinking $Y$ around $W$ suitably, the log canonical ring
$$
R(X, \Delta)=\bigoplus_{m \geq 0}\pi_*\mathcal{O}_X(\llcorner m(K_X+\Delta)\lrcorner)
$$
is a locally finitely generated graded $\mathcal{O}_Y$-algebra.
\end{corollary}

\begin{corollary}[{see, \cite[Theorem 4.6]{Fjn12}}] \label{Corollary 1.6}
Let $(X, \Delta)$ be a pair of a normal complex surface $X$ and a boundary $\mathbb{Q}$-divisor $\Delta$ on $X$ such that $K_X+\Delta$ is $\mathbb{Q}$-Cartier.
Let $\pi\colon X\rightarrow Y$ be a projective surjective morphism of complex analytic varieties. 
Assume that one of \textnormal{Cases (C) and (D)} holds.
Then the log canonical ring
$$
R(X, \Delta)=\bigoplus_{m \geq 0}\pi_*\mathcal{O}_X(\llcorner m(K_X+\Delta)\lrcorner)
$$
is a locally finitely generated graded $\mathcal{O}_Y$-algebra.
\end{corollary}

In \cite{Fjn12}, Fujino showed that the analogues of the results in Cases (B), (C), and (D) in the above corollaries are true in the algebraic setting.
In \cite{Hsh16}, Hashizume generalized the finite generation of the log canonical rings. He proved that the finite generation of the adjoint ring for $\mathbb{Q}$-factorial log surfaces also holds.

This paper is organized as follows. 
In Section \ref{section 2}, we recall basic definitions and results, including the openness of nefness.
In Section \ref{section 3}, we prove that the $(K_X+\Delta)$-minimal model program over $Y$ around $W$ works in each of Cases (A), (B), (C), and (D).
In this section, we also treat a case without the assumption that $K_X$ is $\mathbb{Q}$-Cartier. Moreover, in all cases, we prove that the program can be carried out within the corresponding class.
We mainly use Fujino’s cone and contraction theorems (cf.~\cite[Theorems 4.6.1 and 4.6.2]{Fjn24}). 
In Section \ref{section 4}, we prove the abundance theorem when $\Delta$ is a $\mathbb{Q}$-divisor. This section is the main part of this paper. We first treat the case where $K_X+\Delta$ is big over $Y$. In the proof, we apply the minimal model program under the new assumption given in Case (A). Moreover, we use the vanishing theorem (cf.~Lemma \ref{lem: vanishing} below). When $K_X+\Delta$ is not big over $Y$, we reduce the problem to the case when $Y$ is a curve.
In Section \ref{section 5}, we extend the abundance theorem to the case when $\Delta$ is an $\mathbb{R}$-divisor, reducing it to the $\mathbb{Q}$-divisor case by using the technique of Shokurov’s polytope.
In Section \ref{section 6}, we prove the main theorem and corollaries stated in this introduction.

\begin{ack}
The author would like to thank Professor Osamu Fujino for various suggestions, constant support, and warm encouragement.
The author also wishes to thank the referee for their careful reading and for pointing out a gap in the proof of Lemma \ref{Lemma case C Q-Cartier} in the previous version of the manuscript.
\end{ack}

\section{Preliminaries}\label{section 2}
In this section, we will collect some basic definitions and results.
For the details of the definitions of $\mathbb{Q}$-divisors, $\mathbb{R}$-divisors, linear equivalence, Cartier divisors, ampleness, bigness, and pseudo-effectiveness in the complex analytic setting, see \cite{Fjn22} or \cite{Fjn24}.

\begin{definition}
Let $D$ be an $\mathbb{R}$-divisor on a normal complex analytic variety. We say that $D$ is \textit{boundary} if every coefficient $a_i$ of $D=\sum a_iD_i$ satisfies $a_i\in[0, 1]$. 
We put
\begin{align*}
    \llcorner D\lrcorner&:=\sum\llcorner a_i\lrcorner D_i,
    \quad \ulcorner D\urcorner:=\sum\ulcorner a_i\urcorner D_i,
    \quad \{ D\}:=\sum(a_i-\llcorner a_i\lrcorner) D_i,\\
    D^{>1}&:=\sum_{a_i>1}a_i D_i,
    \quad D^{=1}:=\sum_{a_i=1}a_i D_i,\quad\text{and}
    \quad D^{<1}:=\sum_{a_i<1}a_i D_i.
\end{align*}
\end{definition}

\begin{remark}[Stein factorization]\label{rem: Stein factorization}
We can take the Stein factorization of each proper surjective morphism of reduced complex analytic spaces (cf.~\cite[Theorem 1.9]{Uen75}).
\end{remark}

\begin{definition}[{Nefness, see \cite[Definition 2.48]{Fjn22}}]
Let $\pi\colon X\rightarrow Y$ be a projective morphism of complex analytic spaces such that $X$ is a normal complex analytic variety, and let $W$ be a subset of $Y$. Let $D$ be an $\mathbb{R}$-Cartier $\mathbb{R}$-divisor on $X$. If $D\cdot C\geq0$ for every projective integral curve $C$ on $X$ such that $\pi(C)$ is a point, then $D$ is said to be \textit{$\pi$-nef} or \textit{nef over $Y$}. If $D\cdot C\geq0$ for every projective integral curve $C$ on $X$ such that $\pi(C)$ is a point of $W$, then $D$ is said to be \textit{nef over $W$} or \textit{$\pi$-nef over $W$}.
\end{definition}

\begin{definition}[Suitable neighborhoods] \label{defn: suitable neighborhood}
Let $\pi\colon X\rightarrow Y$ be a projective surjective morphism of complex analytic varieties. Assume that $X$ is a normal complex surface. Let $W$ be a compact subset of $Y$.
If an open neighborhood $U$ of $W$ on $Y$ satisfies the following properties, we call it a \textit{suitable neighborhood of $W$ associated to $\pi\colon X\rightarrow Y$.}
\begin{itemize}
    \item $U$ is a relatively compact subset of $Y$.
    \item $U-W$ is smooth.
    \item $\pi\colon X\rightarrow Y$ is smooth over $U-W$.
\end{itemize}
\end{definition}

\begin{remark}
The definition given in Definition \ref{defn: suitable neighborhood} is not standard. We introduced it for convenience.
\end{remark}

\begin{lemma}[Openness of nefness] \label{lem: openness of nef}
Let $X$ be a normal complex surface and let $\pi\colon X\rightarrow Y$ be a projective surjective morphism of complex analytic varieties.
Let $W$ be a compact subset of $Y$.
Let $U$ be a suitable neighborhood of $W$ associated to $\pi\colon X\rightarrow Y$ in the sense of Definition \ref{defn: suitable neighborhood}.
If an $\mathbb{R}$-Cartier $\mathbb{R}$-divisor $D$ is nef over $W$, then it is nef over $U$.
In particular, if an $\mathbb{R}$-Cartier $\mathbb{R}$-divisor $D$ is nef over $W$, then it is nef over some open neighborhood of $W$.
\end{lemma}
\begin{proof}
By Remark \ref{rem: Stein factorization}, we can take the Stein factorization, and then we may assume that $\pi$ has connected fibers.
When $\dim Y=0$, we have $U=W=Y$ and there is nothing to prove.
When $\dim Y=1$, since we have replaced $\pi$ with its Stein factorization, $Y$ is a smooth curve and $\pi$ has smooth connected fibers on $U-W$. Then $\pi$ is flat over $Y$, and $\pi$ has irreducible fibers on $U-W$.
Pick a point $P\in W$. Since $\pi$ is flat, we have $\pi^*Q\cdot D=\pi^*P\cdot D\geq0$ for each $Q\in Y$. Since the fiber $\pi^{-1}(Q)$ is irreducible for each $Q\in U-W$, $D$ is $\pi$-nef over $U-W$. Since $D$ is also $\pi$-nef over $W$, $D$ is $\pi$-nef over $U$.
When $\dim Y=2$, since $\pi$ has connected fibers, $\pi$ is isomorphic over $U-W$. Then $D$ is clearly $\pi$-nef over $U$.
\end{proof}

\begin{remark}
When $\dim X>2$, an $\mathbb{R}$-Cartier $\mathbb{R}$-divisor $D$ which is nef over $W$ is not necessarily nef over a neighborhood of $W$. This fact often causes troublesome problems when we discuss the minimal model program for projective morphisms between complex analytic spaces (cf.~\cite[Remark~1.5]{Fjn22}). However, in the case when $\dim X=2$, there is no such problem by Lemma \ref{lem: openness of nef}.
\end{remark}

\begin{definition}[{Semi-ampleness, see \cite[Definition 2.45]{Fjn22}}] 
Let $\pi\colon X\rightarrow Y$ be a projective morphism of complex analytic spaces. A finite $\mathbb{R}_{>0}$-linear (resp.~$\mathbb{Q}_{>0}$-linear) combination of $\pi$-semi-ample Cartier divisors is called a \textit{$\pi$-semi-ample} $\mathbb{R}$-divisor (resp.~$\mathbb{Q}$-divisor). 
\end{definition}

\begin{definition}[{$\mathbb{Q}$-factoriality, see \cite[Definition 2.38]{Fjn22}}]
Let $X$ be a normal complex analytic variety and let $K$ be a compact subset of $X$. Then $X$ is said to be \textit{$\mathbb{Q}$-factorial at $K$} if every prime divisor defined on an open neighborhood of $K$ is $\mathbb{Q}$-Cartier at any point $x\in K$.
Let $\pi\colon X\rightarrow Y$ be a projective morphism and let $W$ be a compact subset of $Y$. If $X$ is $\mathbb{Q}$-factorial at $\pi^{-1}(W)$, then we usually say that $X$ is \textit{$\mathbb{Q}$-factorial over $W$}.
\end{definition}

\begin{definition}[{Kleiman--Mori cones, see \cite[Section 4]{Fjn22}}]
Let $\pi\colon X\rightarrow Y$ be a projective morphism of complex analytic spaces and let $W$ be a compact subset of $Y$. Let $Z_1(X/Y; W)$ be the free abelian group generated by the projective integral curves $C$ on $X$ such that $\pi(C)$ is a point of $W$. Let $U$ be any open neighborhood of $W$. Then we can consider the following intersection pairing
$$
\cdot \colon \operatorname{Pic}(\pi^{-1}(U))\times Z_1(X/Y; W)\rightarrow\mathbb{Z},
$$
given by $\mathcal{L}\cdot C\in\mathbb{Z}$ for $\mathcal{L}\in\operatorname{Pic}(\pi^{-1}(U))$ and $C\in Z_1(X/Y; W)$. We put
$$
\widetilde{A}(U, W):=\operatorname{Pic}(\pi^{-1}(U))/{\equiv _W}
$$
and define
$$
A^1(X/Y; W):=\varinjlim_{W\subset U}\widetilde{A}(U, W),
$$
where $U$ runs through all the open neighborhoods of $W$.
We assume that $A^1(X/Y; W)$ is a finitely generated abelian group. Then we can define the \textit{relative Picard number} $\rho(X/Y; W)$ to be the rank of $A^1(X/Y; W)$. We put
$$
N^1(X/Y; W):=A^1(X/Y; W)\otimes_\mathbb{Z}\mathbb{R}.
$$
Let $A_1(X/Y; W)$ be the image of
$$
Z_1(X/Y; W)\rightarrow \operatorname{Hom}_\mathbb{Z}(A^1(X/Y; W), \mathbb{Z})
$$
given by the above intersection pairing. Then we set
$$
N_1(X/Y; W):=A_1(X/Y; W)\otimes_\mathbb{Z}\mathbb{R}.
$$
As usual, we can define the \textit{Kleiman--Mori cone} $\overline{\text{NE}}(X/Y; W)$ of $\pi\colon X\rightarrow Y$ over $W$, that is, $\overline{\text{NE}}(X/Y; W)$ is the closure of the convex cone in $N_1(X/Y; W)$ spanned by the projective integral curves $C$ on $X$ such that $\pi(C)$ is a point of $W$.
\end{definition}

\begin{lemma}[{Negativity lemma, see \cite[Theorem 4-6-1]{Mts02}}] \label{lem: negativity}
Let $f\colon V\rightarrow S$ be a proper bimeromorphic morphism from a nonsingular complex surface $V$ to a normal complex surface $S$. Let $P\in S$ be a point. Then the intersection matrix $(E_i\cdot E_j)$ is negative definite, where $E_i$'s are the exceptional divisors mapped into $P$.
\end{lemma}

\begin{lemma}[{Zariski's lemma, see \cite[Chapter III, Lemma 8.2]{BHPV04}}] \label{lem: Zariski}
Let $f\colon X\rightarrow C$ be a proper surjective connected morphism from a smooth complex surface $X$ to a smooth curve $C$. Let $f^*y=X_y=\sum_i n_iC_i$, where $n_i>0$ and $C_i\subset X$ are irreducible, be a fiber of $f$. Then we have
\begin{enumerate}
    \item[\textnormal{(1)}] $C_i\cdot X_y=0$ for each $i$,
    \item[\textnormal{(2)}] If $D=\sum_im_iC_i$, where $m_i\in\mathbb{Z}$, then $D^2\leq0$, and
    \item[\textnormal{(3)}] $D^2=0$ holds in \textnormal{(2)} if and only if $D=rX_y$, where $r\in\mathbb{Q}$.
\end{enumerate}
\end{lemma}

\begin{definition}[{Singularities of pairs, see \cite[2.1.1]{Fjn24}}]
We consider a normal complex analytic variety $X$. Let $X_{sm}$ denote the smooth locus of $X$. Let $\Delta$ be an $\mathbb{R}$-divisor on $X$ such that $K_X+\Delta$ is $\mathbb{R}$-Cartier and let $f\colon Y\rightarrow X$ be a proper bimeromorphic morphism from a normal complex analytic variety $Y$. We take a small Stein open subset $U$ of $X$ where $K_U+\Delta|_U$ is a well-defined $\mathbb{R}$-Cartier $\mathbb{R}$-divisor on $U$. In this situation, we can define $K_{f^{-1}U}$ and $K_U$  such that $f_*K_{f^{-1}U}=K_U$. Then we can write
$$
K_{f^{-1}U}=f^*(K_U+\Delta|_U)+E_U
$$
as usual. Note that $E_U$ is a well-defined $\mathbb{R}$-divisor on $f^{-1}U$ such that $f_*E_U=\Delta|_U$. Then we have the following formula
$$
K_Y=f^*(K_X+\Delta)+\sum_E a(E, X, \Delta)E
$$
as in the algebraic case. We note that $\sum_E a(E,X, \Delta)E$ is a globally well-defined $\mathbb{R}$-divisor on $Y$ such that
$$
\left.\left(\sum_E a(E,X, \Delta)E\right)\right|_{f^{-1}U}=E_U
$$
although $K_X$ and $K_Y$ are well-defined only locally.

If $\Delta$ is a boundary $\mathbb{R}$-divisor and $a(E, X, \Delta)\geq-1$ holds for any $f\colon Y\rightarrow X$ and every $f$-exceptional divisor $E$, then $(X, \Delta)$ is called a \textit{log canonical} pair. If $(X,\Delta)$ is log canonical, $\llcorner\Delta\lrcorner=0$, and $a(E, X, \Delta)>-1$ for any $f\colon Y\rightarrow X$ and every $f$-exceptional divisor $E$, then $(X, \Delta)$ is called a \textit{kawamata log terminal} pair.

Let $X$ be a normal complex analytic variety and let $\Delta$ be an effective $\mathbb{R}$-divisor on $X$ such that $K_X+\Delta$ is $\mathbb{R}$-Cartier. The image of $E$ with $a(E, X, \Delta)=-1$ for some $f\colon Y\rightarrow X$ such that $(X, \Delta)$ is log canonical around general points of $f(E)$ is called a \textit{log canonical center} (for short, \textit{lc center}) of $(X, \Delta)$. The image of $E$ with $a(E, X, \Delta)\leq-1$ for some $f\colon Y\rightarrow X$ is called a \textit{non-klt center} of $(X, \Delta)$. The \textit{non-lc locus} of $(X, \Delta)$, denoted by $\operatorname{Nlc}(X, \Delta)$ is the smallest closed subset $Z$ of $X$ such that the complement $(X-Z, \Delta|_{X-Z})$ is log canonical. 
\end{definition}

\begin{lemma}[{Vanishing theorem of Reid--Fukuda type, see \cite[Theorem 1.2]{Fjn23}}] \label{lem: vanishing}
Let $(\widetilde{X}, \Delta)$ be a pair of smooth complex analytic variety $\widetilde{X}$ and a boundary $\mathbb{R}$-divisor $\Delta$ on $\widetilde{X}$ such that $\operatorname{Supp}\Delta$ is a simple normal crossing divisor on $\widetilde{X}$. Let $f\colon \widetilde{X}\rightarrow X$ and $\pi\colon X\rightarrow Y$ be projective morphisms between complex analytic varieties, and let $\mathcal{L}$ be a line bundle on $\widetilde{X}$. If $\mathcal{L}-(\omega_{\widetilde{X}}+\Delta)\sim_{\mathbb{R}}f^*\mathcal{H}$ holds for some $\mathbb{R}$-line bundle $\mathcal{H}$ on $X$ which is nef and log big over $Y$ with respect to $f\colon(\widetilde{X}, \Delta)\rightarrow X$, then
$$
\mathcal{R}^p\pi_*\mathcal{R}^qf_*\mathcal{L}=0
$$
holds for every $p>0$ and every $q$.
We note that an $\mathbb{R}$-line bundle $\mathcal{H}$ on $X$ is called log big over $Y$ with respect to $f\colon(\widetilde{X}, \Delta)\rightarrow X$ if $\mathcal{H}$ is $\pi$-big and $\mathcal{H}|_{f(S)}$ can be written as a finite $\mathbb{R}_{>0}$-linear combination of $\pi$-big line bundles on $f(S)$ for every log canonical center $S$ of $(\widetilde{X}, \Delta)$.
\end{lemma}
See also \cite{FM25}, which discusses a vanishing theorem for surfaces formulated in a different way.

\section{Minimal model program} \label{section 3}
In this section, we will prove the minimal model program.
\begin{lemma}[Cone and contraction theorem] \label{Lemma 3.1}
Let $\pi\colon X\rightarrow Y$ be a projective surjective morphism of complex analytic varieties such that $X$ is a normal complex surface, and let $W$ be a compact subset of $Y$ such that $\rho(X/Y; W)<\infty$. Let $\Delta$ be an effective $\mathbb{R}$-divisor on $X$ such that $K_X+\Delta$ is $\mathbb{R}$-Cartier.
\begin{enumerate}
\item[\textnormal{(1)}] $\overline{\textnormal{NE}}(X/Y; W)=\overline{\textnormal{NE}}(X/Y; W)_{(K_X+\Delta)\geq0}+\sum R_j$, where $R_j$'s are rational $(K_X+\Delta)$-negative extremal rays of $\overline{\textnormal{NE}}(X/Y; W)$.
\item[\textnormal{(2)}] For each rational $(K_X+\Delta)$-negative extremal face $F$ of\\ $\overline{\textnormal{NE}}(X/Y; W)$, after shrinking $Y$ around $W$ suitably, there is a projective morphism $\varphi_F\colon X\rightarrow Z$ over $Y$ with the following properties.
\begin{enumerate}
\item[\textnormal{(i)}] Let $C$ be a projective curve on $X$ such that $\pi(C)$ is a point of $W$. Then $\varphi_F(C)$ is a point if and only if the numerical equivalence class $[C]$ of $C$ is in $F$.
\item[\textnormal{(ii)}] The natural map $\mathcal{O}_Z\rightarrow(\varphi_F)_*\mathcal{O}_X$ is an isomorphism.
\item[\textnormal{(iii)}] Let $\mathcal{L}$ be a line bundle on $X$ such that $\mathcal{L}\cdot C=0$ for every curve $C$ with $[C]\in F$. After shrinking $Y$ around $W$ again, there exists a line bundle $\mathcal{L}_Z$ on $Z$ such that $\mathcal{L}\cong\varphi_F^*\mathcal{L}_Z$ holds. 
\end{enumerate}
\end{enumerate}
\end{lemma}
\begin{proof}
This is obvious by \cite[Theorem 4.6.1, Theorem 4.6.2]{Fjn24}. 
Since $\operatorname{dim}\operatorname{Nlc}(X, \Delta)=0$, 
we have $\overline{\text{NE}}(X/Y; W)_{\operatorname{Nlc}(X, \Delta)}=0$. 
\end{proof}

We note a property of $(K_X+\Delta)$-negative extremal curves in Case (C).
\begin{lemma} \label{lem: Q-Cartier}
Under the assumptions of Lemma \ref{Lemma 3.1}, we further assume that \textnormal{Case (C)} holds.
Let $C$ be an irreducible projective curve which is mapped into a point on $W$. Suppose that there exists the contraction morphism $\varphi_R\colon X\rightarrow Z$ associated to a rational $(K_X+\Delta)$-negative extremal ray $R$ of $\operatorname{\overline{NE}(X/Y; W)}$,  given by Lemma \ref{Lemma 3.1}, such that
$\varphi_R(C)$ is a point.
Suppose that $\operatorname{dim}Y\geq1$.
Then, $C$ is $\mathbb{Q}$-Cartier.
\end{lemma}
\begin{proof}
By Lemma \ref{Lemma 3.1}, 
we have $[C]\in R$.
Since $\operatorname{dim}Y\geq1$, we have $\operatorname{dim}Z\geq1$.
When $\operatorname{dim}Z=1$, then $Z$ is a smooth curve.
If $\operatorname{Supp}\varphi_R^*(\varphi_R(C))$ is irreducible, $C$ is proportional to the pullback of a point in $Z$.
In particular, $C$ is $\mathbb{Q}$-Cartier.
Suppose now that $\operatorname{Supp}\varphi_R^*(\varphi_R(C))$ is reducible. (Actually, it turns out that this situation does not occur. See the proof of the next Lemma.)
In this case, we have $C^2<0$ in the sense of Mumford (see, for example, \cite[Section 1]{Fum84}), because $\varphi_R$ has connected fibers.
On the other hand, when $\operatorname{dim}Z=2$, we also have $C^2<0$ in the sense of Mumford.
Then there exists a bimeromorphic contraction morphism $\psi\colon X\rightarrow X'$ such that $\operatorname{Exc}(\psi)=C$, by the Grauert--Sakai contraction theorem (cf.~\cite[Theorem 1.2]{Fum84}).
By Lemma \ref{lem: negativity}, $(X', \psi_*\Delta)$ is numerically log canonical (for the definition of numerical log canonicity, see \cite[Notation 4.1]{KM98}). This implies that $K_{X'}+\psi_*\Delta$ is $\mathbb{R}$-Cartier.
Thus we can write $K_X+\Delta=\psi^*(K_{X'}+\psi_*\Delta)+aC$, where $a\in\mathbb{R}_{>0}$.
In particular, $C$ is $\mathbb{Q}$-Cartier.
\end{proof}

We now take a closer look at the contraction associated with the extremal ray.
\begin{lemma} \label{Lemma case C Q-Cartier}
Under the assumptions of Lemma \ref{Lemma 3.1}, we further assume that one of \textnormal{Cases (A)-(D)} holds.
We replace $Y$ with a suitable neighborhood $U$ of $W$ associated to $\pi\colon X\rightarrow Y$ in the sense of Definition \ref{defn: suitable neighborhood}.
Then, we need not shrink $Y$ anymore in Lemma \ref{Lemma 3.1}, when the extremal face $F$ in Lemma \ref{Lemma 3.1} is a ray.
More precisely, for each rational $(K_X+\Delta)$-negative extremal ray $R$ of $\overline{\textnormal{NE}}(X/Y; W)$, without shrinking $Y$ around $W$ anymore, there is a projective morphism $\varphi_R\colon X\rightarrow Z$ over $Y$ with the following properties.
\begin{enumerate}
\item[\textnormal{(i)}] Let $C$ be a projective curve on $X$ such that $\pi(C)$ is a point of $W$. Then $\varphi_R(C)$ is a point if and only if the numerical equivalence class $[C]$ of $C$ is in $R$.
\item[\textnormal{(ii)}] The natural map $\mathcal{O}_Z\rightarrow(\varphi_R)_*\mathcal{O}_X$ is an isomorphism.
\item[\textnormal{(iii)'}] Let $\mathcal{L}$ be a line bundle on $X$ such that $\mathcal{L}\cdot C=0$ for every curve $C$ with $[C]\in R$. Without shrinking $Y$ around $W$ anymore, there exists a line bundle $\mathcal{L}_Z$ on $Z$ such that $\mathcal{L}\cong\varphi_R^*\mathcal{L}_Z$ holds. 
\item[\textnormal{(iv)}] $Y$ is a suitable neighborhood of $W$ also associated to the induced morphism $\pi_Z\colon Z \rightarrow Y$ in the sense of Definition \ref{defn: suitable neighborhood}.
\end{enumerate}
\end{lemma}

\begin{proof}
By Lemma \ref{Lemma 3.1}, we have a contraction morphism $\phi_R$ over some open neighborhood of $W$.
We divide the argument according to the dimensions of $Y$ and $\operatorname{Im} \phi_R$.

\case
When $\operatorname{dim}Y=0$, we have $U=W=Y$. Then this lemma is obvious by Lemma \ref{Lemma 3.1}.   

\medskip
We note that when $\operatorname{dim}Y\geq1$, every exceptional curve of $\phi_R$ is $\mathbb{Q}$-Cartier. Indeed, it is clear in Cases (A) and (B), and follows from Lemma \ref{lem: Q-Cartier} in Case (C). In Case (D), any divisor on an open subset of $X$ is $\mathbb{Q}$-Cartier  (cf.~\cite[Theorem 7.3.2]{Ish18}).
\par
In Cases \ref{Case 2} and \ref{Case 3}, we assume that $\dim Y=1$. In these cases, taking the Stein factorization of $\pi$, we may assume that the fiber $\pi^{-1}(P)$ is smooth and irreducible for each $P\in Y-W$. 
Let $C$ be a curve on $X$ mapped into a point on $W$ such that the numerical equivalence class $[C]\in\overline{\textnormal{NE}}(X/Y; W)$ generates $R$.

\case \label{Case 2}
When $\dim Y=1$ and $\dim (\operatorname{Im} \phi_R)=1$, we can write $\phi_R^*(\phi_R(C))=E=\sum a_iE_i$, where $a_i>0$ and $E_i$'s are irreducible curves.
As noted above, each $E_i$ is a $\mathbb{Q}$-Cartier divisor. If $\operatorname{Supp}E\neq\operatorname{Supp}C$, then some $E_i$ satisfies $E_i\cdot C>0$ because $\phi_R$ has connected fibers. By the property (i) in Lemma \ref{Lemma 3.1}, we have $E^2>0$. This is a contradiction, by Lemma \ref{lem: Zariski}.
Then, $\operatorname{Supp}\phi_R^*(\phi_R(C))=\operatorname{Supp}C$.
By Lemma \ref{lem: Zariski}, we have $C^2=0$. Then, $C$ is a support function of $R$. Moreover, since $C$ is proportional to a pullback of a point in $Z$, $mC$ is $\pi$-free over $Y$ for some positive integer $m$.
We obtain a morphism $X\rightarrow\mathbb{P}_Y(\pi_*\mathcal{O}_X(mC))$. Taking the Stein factorization of this morphism, we obtain a contraction morphism $\varphi_R\colon X\rightarrow Z$ over $Y$. This $\varphi_R$ coincides with $\phi_R$ on the open subset of $X$ on which $\phi_R$ is defined.
This morphism satisfies the properties (i), (ii), and (iv) of this lemma.

To show that $\varphi_R$ satisfies (iii)', take a line bundle $\mathcal{L}$ on $X$ such that $\mathcal{L}\cdot C=0$. 
We take a general fiber $F$ of $\varphi_R$. We may assume that $F$ is smooth, by Sard's theorem (cf.~\cite[Corollary 1.8]{Uen75}). Since $-(K_X+\Delta)$ is $\varphi_R$-ample, $-(K_X+\Delta)|_F=-K_F-\Delta|_F$ is ample. Then, $-K_F$ is ample. This implies that the genus of the smooth curve $F$ is 0, that is, $F\cong\mathbb{P}^1$.
Since $\varphi_R$ is flat, $\mathcal{L}\cdot F=0$. So, $\mathcal{L}|_F\cong\mathcal{O}_{\mathbb{P}^1}$.
Then, for a general point $z\in Z$, we have
\begin{align*}
    h^0(X_z, \mathcal{L}|_{X_z})=1,\quad\text{and}\quad
    h^0(X_z, \mathcal{L}^{-1}|_{X_z})=1.
\end{align*}
By the upper semi-continuity (cf.~\cite[Theorem 1.4]{Uen75}), for every point $z\in Z$, we have
\begin{align*}
    h^0(X_z, \mathcal{L}|_{X_z})\geq1,\quad\text{and}\quad
    h^0(X_z, \mathcal{L}^{-1}|_{X_z})\geq1.
\end{align*}
This implies that $\mathcal{L}|_{X_z}\cong\mathcal{O}_{X_z}$ for every point $z\in Z$. 
In particular, $\mathcal{L}$ is $\varphi_R$-free over $Z$.
Since $\mathcal{L}$ is $\varphi_R$-numerically trivial, it is a pull-back of a line bundle on $Z$. Then, $\varphi_R$ satisfies (iii)'.

\case \label{Case 3}
When $\dim Y=1$ and $\dim (\operatorname{Im}\phi_R)=2$, we have $C^2<0$ by Lemma \ref{lem: negativity}. 
Take a $\pi$-ample line bundle $\mathcal{A}$ on $X$. Then there is a positive rational number $\lambda$ such that $C+\lambda \mathcal{A}$ is a support function of $R$. 
Put $H:=C+\lambda \mathcal{A}$.
We can take some positive integer $a$ such that $aH-(K_X+\Delta)$ is $\pi$-ample over $W$. By the basepoint-free theorem (cf.~\cite[Theorem 1.1.11]{Fjn24}), $mH$ is $\pi$-free over a neighborhood of $W$ for every sufficiently large positive integer $m$.
By the definition of $H$, for a point $P\in Y-W$, we have $H\cdot\pi^*P>0$. Then $m_1(H\cdot \pi^*P)\geq 2g+1$ for every sufficiently large positive integer $m_1$, where $g$ is the genus of the smooth curve $\pi^{-1}(P)$. Then $m_1H|_{\pi^{-1}(P)}$ is very ample. Since $\pi\colon X\rightarrow Y$ is flat, $H\cdot\pi^*Q=H\cdot\pi^*P$ and the genus of $\pi^{-1}(Q)$ is equal to the genus $g$ of $\pi^{-1}(P)$ for every point $Q\in Y-W$. Then $m_1H|_{\pi^{-1}(Q)}$ is very ample for every point $Q\in Y-W$. This implies that $mH$ is $\pi$-very ample over $Y-W$ for every sufficiently large positive integer $m$.
Then, $mH$ is $\pi$-free over $Y$ for a positive integer $m$, and we obtain a morphism $\varphi_R\colon X\rightarrow Z$ over $Y$ as in Case \ref{Case 2}. This $\varphi_R$ coincides with $\phi_R$ on the open subset of $X$ on which $\phi_R$ is defined. This morphism satisfies the properties (i), (ii), and (iv) of this lemma.

To show that $\varphi_R$ satisfies (iii)', we take a line bundle $\mathcal{L}$ on $X$ such that $\mathcal{L}\cdot C=0$. 
Then $\mathcal{L}-(K_X+\Delta)$ is $\varphi_R$-ample over $\pi_Z^{-1}W$. By the basepoint-free theorem (cf.~\cite[Theorem 1.1.11]{Fjn24}), $m\mathcal{L}$ is $\varphi_R$-free over a neighborhood of $\pi_Z^{-1}W$ for every sufficiently large positive integer $m$.

We will check that $m\mathcal{L}$ is $\varphi_R$-free over $Z-\pi_Z^{-1}W$ for every sufficiently large positive integer $m$.
Since $\dim Z=2$, let $E$ be a divisor on $X$ such that $\operatorname{Supp}E$ corresponds to the exceptional locus of $\varphi_R$. By the property (i), the numerical equivalence class $[E_i]\in\overline{\textnormal{NE}}(X/Y; W)$ of each irreducible component $E_i$ of $E$ is included in $R$. Since $E_i\cdot H=0$, we have $E_i\subset \operatorname{Supp}C$. Then, the exceptional locus of $\varphi_R$ corresponds to $\operatorname{Supp}C$.
Since $\varphi_R$ is isomorphic over $Z-\pi_Z^{-1}W$, $m\mathcal{L}$ is $\varphi_R$-free over $Z-\pi_Z^{-1}W$ for any positive integer $m$.

We have checked that $m\mathcal{L}$ is $\varphi_R$-free over $Z$ for every sufficiently large positive integer $m$. 
Then, there is some positive integer $m_0$ such that both $m_0\mathcal{L}$ and $(m_0+1)\mathcal{L}$ are $\varphi_R$-free over $Z$. Since both $m_0\mathcal{L}$ and $(m_0+1)\mathcal{L}$ are $\varphi_R$-numerically trivial, they are pull-backs of line bundles on $Z$. Therefore, $\mathcal{L}$ is also a pull-back of a line bundle on $Z$. Then, $\varphi_R$ satisfies (iii)'.

\case
When $\dim Y=2$, taking the Stein factorization of $\pi$, we may assume that $\pi\colon X\rightarrow Y$ is isomorphic over $Y-W$. Then any line bundle on $X$ is $\pi$-free over $Y-W$. We can show this case in the same way as in Case \ref{Case 3}.
\end{proof}

\begin{lemma} \label{lem: preserve conditions}
Assume that $(X, Y, W, \pi, \Delta)$ satisfies the condition \textnormal{($\bigstar$)} and that one of \textnormal{Cases (A)-(D)} holds.
Let $\varphi\colon X\rightarrow X'$ be a bimeromorphic contraction morphism associated to a $(K_X+\Delta)$-negative ray of $\overline{\textnormal{NE}}(X/Y; W)$, which is given by Lemma \ref{Lemma case C Q-Cartier}. Let $\pi'\colon X'\rightarrow Y$ be the induced morphism and let $\Delta':=\varphi_*\Delta$. 

Then $(X', Y, W, \pi', \Delta')$ also satisfies the condition \textnormal{($\bigstar$)}. Moreover, $K_{X'}$ is $\mathbb{Q}$-Cartier and every prime divisor on $X'$ which is mapped into $W$ is $\mathbb{Q}$-Cartier in Case \textnormal{(A)}, $X'$ is $\mathbb{Q}$-factorial over $W$ in Case \textnormal{(B)}, $(X', \Delta')$ is log canonical in Case \textnormal{(C)}, and $X'$ has only rational singularities in Case \textnormal{(D)}.
\end{lemma}
\begin{proof}
We note that the contraction morphism $\varphi$ preserves normality by (ii) of Lemma \ref{Lemma case C Q-Cartier}.
We also note that $\rho(X/Y; W)-\rho(X'/Y; W)=1$, because we have the following exact sequence
$$
0\rightarrow A^1(X'/Y; W)\xrightarrow{D\mapsto\varphi^*D}
A^1(X/Y; W)\xrightarrow{D\mapsto(D\cdot C)}
\mathbb{Z},
$$
where $C$ is a curve on $X$ such that its numerical equivalence class $[C]\in\overline{\text{NE}}(X/Y; W)$ generates the extremal ray associated to $\varphi$, by (iii)' of Lemma \ref{Lemma case C Q-Cartier}.

In order to show that $(X', Y, W, \pi', \Delta')$ satisfies the condition ($\bigstar$), we only have to show that $K_{X'}+\Delta'$ is $\mathbb{R}$-Cartier because $(X, Y, W, \pi, \Delta)$ satisfies the condition ($\bigstar$) and $\varphi\colon X\rightarrow X'$ is a bimeromorphic contraction morphism. 

In Case (A), let $E$ be a divisor on $X$ such that its support $\operatorname{Supp}E$ corresponds to the exceptional locus of $\varphi$. Let $E_i$ be an irreducible component of $E$.
By assumption, each $E_i$ is $\mathbb{Q}$-Cartier. Since the contraction morphism $\varphi$ satisfies the conditions in Lemma \ref{Lemma case C Q-Cartier}, the numerical equivalence class $[E_i]\in\overline{\textnormal{NE}}(X/Y; W)$ is proportional to $[C]\in\overline{\textnormal{NE}}(X/Y; W)$. By assumption, $C$ is also $\mathbb{Q}$-Cartier. Since $C^2<0$ implies $C\cdot E_i<0$, we have $\operatorname{Supp}E_i=\operatorname{Supp}C$. Then $\operatorname{Supp}E=\operatorname{Supp}C$ and in particular the exceptional locus of $\varphi$ is irreducible. 
Pick a $\mathbb{Q}$-Cartier $\mathbb{Q}$-divisor $\Gamma$ on $X$. Since $\rho(X/Y; W)-\rho(X'/Y; W)=1$ and $E\cdot E<0$, there is a rational number $\lambda$ such that $(\Gamma+\lambda E)\cdot E=0$. Then $\Gamma+\lambda E$ is both $\mathbb{Q}$-Cartier and $\varphi$-numerically trivial. This implies that $\varphi_*(\Gamma+\lambda E)=\varphi_*\Gamma$ is $\mathbb{Q}$-Cartier, by (iii)' of Lemma \ref{Lemma case C Q-Cartier}. 
Since $K_X+\Delta$ is $\mathbb{R}$-Cartier, we can write $K_X+\Delta=\sum_{i=1}^na_i(K_X+\Delta_i)$ on some open neighborhood $V$ of $\operatorname{Supp}E$ on $X$, where every $\Delta_i$ is a boundary $\mathbb{Q}$-divisor such that $K_X+\Delta_i$ is $\mathbb{Q}$-Cartier and every $a_i$ is a positive real number such that $\sum a_i=1$. Since $\operatorname{Supp}E$ corresponds to the exceptional locus of $\varphi$, $\varphi(V)\subset X'$ is open and $V=\varphi^{-1}(\varphi(V))$. Substituting $\Gamma=K_X+\Delta_i$, we see that $\varphi_*(K_X+\Delta_i)$ is $\mathbb{Q}$-Cartier for every $i$. Then, $K_{X'}+\Delta'=\sum_{i=1}^na_i(K_{X'}+\varphi_*\Delta_i)$ is $\mathbb{R}$-Cartier on $\varphi(V)$. Since $\varphi$ is isomorphic on $X-\operatorname{Supp}E$, $K_{X'}+\Delta'$ is $\mathbb{R}$-Cartier on $X'$. Substituting $\Gamma=K_X$, $K_{X'}$ is $\mathbb{Q}$-Cartier.
Pick a prime divisor $D'$ on $X'$ which is mapped into $W$. Let $D$ be the strict transform of $D'$. 
Then $\pi(D)\subset W$ implies that $D$ is $\mathbb{Q}$-Cartier. Then, substituting $\Gamma=D$, $\varphi_*D=D'$ is $\mathbb{Q}$-Cartier.  
Therefore, $(X', Y, W, \pi', \Delta')$ satisfies the condition ($\bigstar$) and $K_{X'}$ is $\mathbb{Q}$-Cartier and every prime divisor on $X'$ which is mapped into $W$ is $\mathbb{Q}$-Cartier.

In Case (B), as above, the exceptional locus of $\varphi$ is irreducible. In order to prove that $X'$ is $\mathbb{Q}$-factorial over $W$, pick a prime divisor $D'$ defined on an open neighborhood $U$ of $(\pi')^{-1}(W)$. Let $D$ be the strict transform of $D'$. Then $D$ is $\mathbb{Q}$-Cartier over a neighborhood of $W$ by the assumption of (B). Then as above, $\varphi_*D=D'$ is $\mathbb{Q}$-Cartier over a neighborhood of $W$. Then $X'$ is $\mathbb{Q}$-factorial over $W$. We can show that $K_{X'}+\Delta'$ is $\mathbb{R}$-Cartier in the same way as above. Therefore, $(X', Y, W, \pi', \Delta')$ satisfies the condition ($\bigstar$) and $X'$ is $\mathbb{Q}$-factorial over $W$.

In Case (C), by Lemma \ref{lem: negativity}, $(X', \Delta')$ is numerically log canonical (for the definition of numerical log canonicity, see \cite[Notation 4.1]{KM98}). This implies that $K_{X'}+\Delta'$ is $\mathbb{R}$-Cartier. Therefore, $(X', Y, W, \pi', \Delta')$ satisfies the condition ($\bigstar$) and $(X', \Delta')$ is log canonical.

In Case (D), $X'$ also has only rational singularities by the vanishing theorem (cf.~\cite[Theorem 3.3.1]{Fjn24}), because the same argument as \cite[Proposition 3.7]{Fjn12} is available, and we have
$$
\mathcal{R}^j\varphi_*\mathcal{O}_X=0
$$
for all $j>0$. In particular, any divisor on an open subset of $X$ is $\mathbb{Q}$-Cartier (cf.~\cite[Theorem 7.3.2]{Ish18}). This implies that $K_{X'}+\Delta'$ is $\mathbb{R}$-Cartier. Therefore, $(X', Y, W, \pi', \Delta')$ satisfies the condition ($\bigstar$) and $X'$ also has only rational singularities.
\end{proof}

\begin{remark}[{See also Remark \ref{rem: A' MMP}}] \label{rem: A' preserves conditions}
If $(X, Y, W, \pi, \Delta)$ satisfies the conditions $(\bigstar)$ and
\begin{enumerate}
    \item[(A')] every prime divisor on $X$ which is mapped into $W$ is $\mathbb{Q}$-Cartier,
\end{enumerate}
then $(X', Y, W, \pi', \Delta')$ also satisfies the conditions $(\bigstar)$ and every prime divisor on $X'$ which is mapped into $W$ is $\mathbb{Q}$-Cartier. This can be proved in the same way as in the proof of Case (A).
\end{remark}

\begin{theorem}[{Minimal model program, see \cite[Theorem 3.3]{Fjn12}}] \label{thm: MMP}
Assume that $(X, Y, W, \pi, \Delta)$ satisfies the condition \textnormal{($\bigstar$)} and that one of \textnormal{Cases (A)-(D)} holds.
We shrink $Y$ around $W$ suitably.
Then, there is a sequence of at most $\rho(X/Y;W)-1$ contractions
$$
(X, \Delta)=(X_0, \Delta_0)\stackrel{\varphi_0}{\rightarrow}(X_1, \Delta_1)\stackrel{\varphi_1}{\rightarrow}\cdots\stackrel{\varphi_{k-1}}{\rightarrow}(X_k, \Delta_k)=(X^*, \Delta^*)
$$
over $Y$ such that
\begin{itemize}
\item \textnormal{(Minimal model)} if $K_X+\Delta$ is $\pi$-pseudo-effective, then $K_{X^*}+\Delta^*$ is nef over $W$, and
\item \textnormal{(Mori fiber space)} if $K_X+\Delta$ is not $\pi$-pseudo-effective, then $(X^*, \Delta^*)$ is a Mori fiber space.
\end{itemize}
We note that $K_{X_i}$ is $\mathbb{Q}$-Cartier and every prime divisor on $X_i$ which is mapped into $W$ is $\mathbb{Q}$-Cartier for every $i$ in \textnormal{Case (A)}, that $X_i$ is $\mathbb{Q}$-factorial over $W$ for every $i$ in \textnormal{Case (B)}, that $(X_i, \Delta_i)$ is log canonical for every $i$ in \textnormal{Case (C)}, and that $X_i$ has only rational singularities for every $i$ in \textnormal{Case (D)}.
\end{theorem}
\begin{proof}
We consider the following step-by-step procedure.
\begin{enumerate}[Step \arabic*.]
\item We start with $(X, \Delta)$.
\item If $K_X+\Delta$ is nef over $W$, we stop.
\item If $K_X+\Delta$ is not nef, then there is a rational $(K_X+\Delta)$-negative extremal ray $R$ of $\overline{\text{NE}}(X/Y; W)$ and we have the corresponding contraction morphism $\varphi_R\colon X\rightarrow X'$ over $Y$ in Lemma \ref{Lemma case C Q-Cartier}. 
We note that $X'$ is normal, and that $\rho(X/Y; W)-\rho(X'/Y; W)=1$ as shown in the proof of Lemma \ref{lem: preserve conditions}.
\item If $\varphi_R\colon X\rightarrow X'$ is not bimeromorphic, then it is a Mori fiber space and we stop.
\item Otherwise, $\varphi_R\colon X\rightarrow X'$ is bimeromorphic. Let $\pi'\colon X'\rightarrow Y$ be the induced morphism and let $\Delta':={\varphi_R}_*\Delta$. 
By Lemma \ref{lem: preserve conditions}, $(X', Y, W, \pi', \Delta')$ also satisfies the condition \textnormal{($\bigstar$)}. Moreover, $K_{X'}$ is $\mathbb{Q}$-Cartier and every prime divisor on $X'$ which is mapped into $W$ is $\mathbb{Q}$-Cartier in Case \textnormal{(A)}, $X'$ is $\mathbb{Q}$-factorial over $W$ in Case \textnormal{(B)}, $(X', \Delta')$ is log canonical in Case \textnormal{(C)}, and $X'$ has only rational singularities in Case \textnormal{(D)}.
Then, we can replace $(X, \Delta)$ with $(X', {\varphi_R}_*\Delta)$ and go back to the first procedure.
\end{enumerate}
We may repeat this process. Since $\rho(X_i/Y; W)-\rho(X_{i+1}/Y; W)=1$, we obtain a sequence of at most $\rho(X/Y;W)-1$ contractions
$$
(X, \Delta)=(X_0, \Delta_0)\stackrel{\varphi_0}{\rightarrow}(X_1, \Delta_1)\stackrel{\varphi_1}{\rightarrow}\cdots\stackrel{\varphi_{k-1}}{\rightarrow}(X_k, \Delta_k)=(X^*, \Delta^*)
$$
over $Y$ such that either $K_{X^*}+\Delta^*$ is nef over $W$ or $(X^*, \Delta^*)$ is a Mori fiber space.

By Lemma \ref{lem: negativity}, we can write
$$
K_X+\Delta=\Phi^*(K_{X^*}+\Delta^*)+E,
$$
where $\Phi\colon X\rightarrow X^*$ is the composition of the $\varphi_i$'s, and $E$ is an effective $\Phi$-exceptional divisor on $X$ whose support corresponds to the exceptional locus of $\Phi$.
If $K_{X^*}+\Delta^*$ is nef over $W$, then it is nef over $Y$ by Lemma \ref{lem: openness of nef}.
Then $K_X+\Delta$ is $\pi$-pseudo-effective.
Conversely, if $K_X+\Delta$ is $\pi$-pseudo-effective, then so is $K_{X^*}+\Delta^*$, because $\Phi\colon X\rightarrow X^*$ is bimeromorphic. In order to show that $(X^*, \Delta^*)$ is a minimal model, we suppose that $\varphi_k\colon (X_k, \Delta_k)=(X^*, \Delta^*)\rightarrow Z$ is a Mori fiber space over $Y$. In this case, $-(K_{X_k}+\Delta_k)$ is $\varphi_k$-ample and $\operatorname{dim}Z<2$. For every point $P\in X_k$, there exists a curve $C$ on $X_k$ passing through $P$ which is contracted by $\varphi_k$, which satisfies that $(K_{X_k}+\Delta_k)\cdot C<0$. This contradicts the fact that $(K_{X_k}+\Delta_k)\cdot C\geq0$ for a very general curve $C$ on $X_k$. Therefore, $(X^*, \Delta^*)$ is a minimal model when $K_X+\Delta$ is $\pi$-pseudo-effective and it is a Mori fiber space when $K_X+\Delta$ is not $\pi$-pseudo-effective.
\end{proof}

\begin{remark} \label{rem: A' MMP}
If $(X, Y, W, \pi, \Delta)$ satisfies the conditions $(\bigstar)$ and
\begin{enumerate}
    \item[(A')] every prime divisor on $X$ which is mapped into $W$ is $\mathbb{Q}$-Cartier,
\end{enumerate}
then we can run the $(K_X+\Delta)$-MMP over $Y$ around $W$ as above. Moreover, $(X_i, Y, W, \pi_i, \Delta_i)$ satisfies the condition $(\bigstar)$ and every prime divisor on $X_i$ which is  mapped into $W$ is $\mathbb{Q}$-Cartier for every $i$. This can be proved in the same way as in the proof of Theorem \ref{thm: MMP}, by Remark \ref{rem: A' preserves conditions}.
\end{remark}

\section{\texorpdfstring{Abundance for $\mathbb{Q}$-divisors}{Abundance for Q-divisors}}\label{section 4}
In this section, we will prove the abundance theorem for $\mathbb{Q}$-divisors.
\begin{lemma}[Adjunction] \label{lem: adjunction}
Let $X$ be a normal complex surface, and $A_j$ a reduced irreducible closed subvariety of dimension 1. Then we have the short exact sequence
$$
0 \rightarrow \mathcal{T} \rightarrow \mathcal{O}_X(K_X+{A_j}) \otimes \mathcal{O}_{A_j} \rightarrow \omega_{A_j} \rightarrow 0
$$
where $\mathcal{T}$ is a sheaf whose support has dimension zero.
\end{lemma}
\begin{proof}[Proof of Lemma \ref{lem: adjunction}]
The corresponding statement for algebraic surfaces holds by \cite[Lemma 4.4]{Fjn12}. Using the technique from \cite{RRV71}, the analogous statement in the analytic setting also holds.
\end{proof}

\begin{theorem}[{see, \cite[Theorem 4.1]{Fjn12}}] \label{Theorem 4.2}
Assume that $(X, Y, W, \pi, \Delta)$ satisfies the condition \textnormal{($\bigstar$)}, $\Delta$ is a $\mathbb{Q}$-divisor, and \textnormal{Case (A)} holds.
Assume that $W:=\{P\}$, where $P\in Y$ is a point. If $K_X+\Delta$ is nef over $W$ and big over $Y$, then it is semi-ample over a neighborhood of $W$.  
\end{theorem}
\begin{remark}
The assumptions in Theorem \ref{Theorem 4.2} are slightly different from \cite[Theorem 4.1]{Fjn12}. Let us explain this difference and the motivation of the assumptions in Theorem \ref{Theorem 4.2}.
The corresponding statement in the algebraic case was proved under the assumption that $X$ is $\mathbb{Q}$-factorial (cf.~\cite[Theorem 4.1]{Fjn12}). 
In the analytic case, for a compact subset $W'$ of $Y$ with $W'\subset W$, $X$ is not necessarily $\mathbb{Q}$-factorial over $W'$ even if $X$ is $\mathbb{Q}$-factorial over $W$. Since semi-ampleness is local on $Y$, this fact causes some problems. In order to avoid such problems, we work under an assumption weaker than $\mathbb{Q}$-factoriality.
\end{remark}

We closely follow Fujino's proof of \cite[Theorem 4.1]{Fjn12} with appropriate modifications. In particular, we adjust the argument to account for the non-compact case.
\begin{proof}[Proof of Theorem \ref{Theorem 4.2}]
By Lemma \ref{lem: openness of nef}, we may assume that $K_X+\Delta$ is nef over $Y$.

\setcounter{stepcounter}{0}
\step \label{step: connected fibers}
We may assume that $\pi$ has connected fibers.    

\begin{proof}[Proof of Step \ref{step: connected fibers}]
By Remark \ref{rem: Stein factorization}, we can take the Stein factorization $X\xrightarrow{\pi'}Y'\xrightarrow{\nu}Y$ of $\pi\colon X\rightarrow Y$. Since $\nu$ is finite, $K_X+\Delta$ is $\pi'$-nef over a neighborhood of $\nu^{-1}W$ and
$$
\nu^*\nu_*\pi'_*\mathcal{O}_X(m(K_X+\Delta))\rightarrow\pi'_*\mathcal{O}_X(m(K_X+\Delta))
$$
is surjective for any positive divisible integer $m$. Since $\pi'^*$ is a right exact functor, 
\begin{align*}
\pi^*\pi_*\mathcal{O}_X(m(K_X+\Delta))&=\pi'^*\nu^*\nu_*\pi'_*{\mathcal{O}_X(m(K_X+\Delta))}\\
&\rightarrow\pi'^*\pi'_*{\mathcal{O}_X(m(K_X+\Delta))}  
\end{align*}
is also surjective. Then if $K_X+\Delta$ is $\pi'$-semi-ample over a neighborhood of $\nu^{-1}W$, it is $\pi$-semi-ample over a neighborhood of $W$. Therefore, we can replace $\pi$ by $\pi'$ and so we may assume that $\pi$ has connected fibers.   
\end{proof}

\step \label{step: separate into A and B}
Taking a relatively compact open neighborhood of $W$, we may assume that $\operatorname{Supp}\Delta$ has only finitely many irreducible components. 
Let $\llcorner\Delta\lrcorner=\sum_iC_i$ be the irreducible decomposition. 
Shrinking $Y$ around $W$, we may assume that if $C_i$ is mapped into a point in $Y$, then it is mapped into $W$.
We put
\begin{align*}
A&:=\sum_{C_i\cdot(K_X+\Delta)=0\text{ and }\pi(C_i)\text{ is a point}}C_i,\textnormal{ and}\\
B&:=\sum_{C_i\cdot(K_X+\Delta)>0\text{ or }\pi(C_i)\text{ is a curve}}C_i.   
\end{align*}
Then $\llcorner\Delta\lrcorner=A+B$.
If $C_i$ is a component of $A$, then we may assume that $C_i^2<0$.

\begin{proof}[Proof of Step \ref{step: separate into A and B}]
When $\operatorname{dim}Y=2$, the dimension of the fibers of $\pi$ is 0 over a Zariski open subset of $Y$. Since $\pi$ has connected fibers, $\pi$ is bimeromorphic. Let $C_i$ be a component of $A$. It suffices to show that the pull-back $C_i'$ of $C_i$ to the resolution of $X$ satisfies $C_i'^2<0$, and so we may assume that $X$ is smooth. Then $C_i^2<0$ by Lemma \ref{lem: negativity}.

When $\operatorname{dim}Y=1$, we may assume that $X$ is smooth as above. Then it is obvious that $C_i^2\leq0$ by Lemma \ref{lem: Zariski}. If $C_i^2=0$, then $C_i$ is proportional to a pull-back of a point in $Y$ by $\pi$, by Lemma \ref{lem: Zariski}. Then $C_i$ is $\pi$-semi-ample and it suffices to show that $K_X+\Delta-C_i$ is semi-ample over a neighborhood of $W$. Replacing $K_X+\Delta$ by $K_X+\Delta-C_i$, the number of components $C_i$ with $C_i^2=0$ decreases.
Therefore, we may assume that $C_i^2<0$ for every component $C_i$.
\end{proof}

From Step \ref{step: 0-dim lc centers} to Step \ref{step: irr comp of B}, we will check that there exists a positive integer $m$ such that the subset
$$
T:=\operatorname{Supp}\operatorname{Coker}(\pi^*\pi_*\mathcal{O}_X(m(K_X+\Delta))\rightarrow\mathcal{O}_X(m(K_X+\Delta)))\subset X
$$
contains no non-klt centers of $(X, \Delta)$, after shrinking $Y$ around $W$ suitably.

\step \label{step: 0-dim lc centers}
Let $Q$ be a zero-dimensional non-klt center of $(X, \Delta)$ such that $Q\notin \operatorname{Supp}A$. Then $Q \notin T$.

\begin{proof}[Proof of Step \ref{step: 0-dim lc centers}]
Let $f\colon \widetilde{X}\rightarrow X$ be a resolution such that $K_{\widetilde{X}}+\Delta_{\widetilde{X}}=f^*(K_X+\Delta)$. We may assume that
\begin{enumerate}
    \item[(1)] $f^{-1}(A)$ has simple normal crossing support, 
    \item[(2)]  $\operatorname{Supp} f_*^{-1} \Delta \cup \operatorname{Exc}(f)$ is a simple normal crossing divisor on $\widetilde{X}$, and
    \item[(3)] $f^{-1}(Q)$ has simple normal crossing support.
\end{enumerate}
Let $W_1$ be the union of the irreducible components of $\Delta_{\widetilde{X}}^{=1}$ which are mapped into $A \cup Q$ by $f$. We put $\Delta_{\widetilde{X}}^{=1}=W_1+W_2$. Then
$$
-W_1-\llcorner\Delta_{\widetilde{X}}^{>1}\lrcorner+\ulcorner-(\Delta_{\widetilde{X}}^{<1})\urcorner-(K_{\widetilde{X}}+\{\Delta_{\widetilde{X}}\}+W_2) \sim_\mathbb{Q}-f^*(K_X+\Delta).
$$
We put
$$
\mathcal{J}_1=f_* \mathcal{O}_{\widetilde{X}}(-W_1-\llcorner\Delta_{\widetilde{X}}^{>1}\lrcorner+\ulcorner-(\Delta_{\widetilde{X}}^{<1})\urcorner)\subset \mathcal{O}_X.
$$
Then we can easily check that
$$
0 \rightarrow \mathcal{J}_1 \rightarrow \mathcal{O}_X(-A) \rightarrow \delta \rightarrow 0
$$
is exact, where $\delta$ is a skyscraper sheaf.
We will check that the conditions in Lemma \ref{lem: vanishing} hold. In other words, we will check that
\begin{itemize}
\item $\mathcal{L}:=\mathcal{O}_{\widetilde{X}}(m(K_{\widetilde{X}}+\Delta_{\widetilde{X}})-W_1-\llcorner\Delta_{\widetilde{X}}^{>1}\lrcorner+\ulcorner-(\Delta_{\widetilde{X}}^{<1})\urcorner)$ is a line bundle on $\widetilde{X}$ for some divisible positive integer $m$,
\item $(\widetilde{X}, \Delta_{\widetilde{X}}-W_1-\llcorner\Delta_{\widetilde{X}}^{>1}\lrcorner+\ulcorner-(\Delta_{\widetilde{X}}^{<1})\urcorner)$ is a pair of smooth complex analytic variety and a boundary $\mathbb{R}$-divisor on it which has a simple normal crossing support,
\item $f\colon \widetilde{X}\rightarrow X$ and $\pi\colon X\rightarrow Y$ are projective, 
\item $\mathcal{L}-\mathcal{O}_{\widetilde{X}}(K_{\widetilde{X}}+\Delta_{\widetilde{X}}-W_1-\llcorner\Delta_{\widetilde{X}}^{>1}\lrcorner+\ulcorner-(\Delta_{\widetilde{X}}^{<1})\urcorner)
=\mathcal{O}_{\widetilde{X}}((m-1)(K_{\widetilde{X}}+\Delta_{\widetilde{X}}))=f^*\mathcal{O}_X((m-1)(K_X+\Delta))$, and
\item $\mathcal{O}_X((m-1)(K_X+\Delta))$ is nef and log big over $Y$ with respect to $f\colon(\widetilde{X}, \Delta_{\widetilde{X}}-W_1-\llcorner\Delta_{\widetilde{X}}^{>1}\lrcorner+\ulcorner-(\Delta_{\widetilde{X}}^{<1})\urcorner)\rightarrow X$.
\end{itemize}
We only have to show the last condition.
We note that an $\mathbb{R}$-line bundle $\mathcal{H}$ on $X$ is log big over $Y$ with respect to $f\colon(\widetilde{X}, \Gamma:=\Delta_{\widetilde{X}}-W_1-\llcorner\Delta_{\widetilde{X}}^{>1}\lrcorner+\ulcorner-(\Delta_{\widetilde{X}}^{<1})\urcorner)\rightarrow X$ if $\mathcal{H}$ is $\pi$-big and $\mathcal{H}|_{f(S)}$ can be written as a finite $\mathbb{R}_{>0}$-linear combination of $\pi$-big line bundles on $f(S)$ for every lc center $S$ of $(\widetilde{X}, \Gamma)$.
By assumption, $K_X+\Delta$ is $\pi$-big.
An lc center $S$ of $(\widetilde{X}, \Gamma)$ is zero-dimensional or a component of $W_2$.
If $f(S)$ is a point, $(K_X+\Delta)|_{f(S)}$ is clearly $\pi$-big.
If $f(S)$ is a curve and it is not contracted by $\pi$, then $\pi|_{f(S)}$ is finite. This implies that $(K_X+\Delta)|_{f(S)}$ is $\pi$-big.
If $f(S)$ is a curve and it is contracted by $\pi$, then $(K_X+\Delta)\cdot f(S)>0$ by the definitions of $B$ and $\Gamma$. This implies that $(K_X+\Delta)|_{f(S)}$ is $\pi$-big.
Therefore, $\mathcal{O}_X((m-1)(K_X+\Delta))$ is log big over $Y$ with respect to $f\colon(\widetilde{X}, \Gamma)\rightarrow X$.
Then we can apply Lemma \ref{lem: vanishing} and we have
$$
\mathcal{R}^p\pi_*(\mathcal{O}_X(m(K_X+\Delta)) \otimes \mathcal{J}_1)=0
$$
for every $p>0$ by Lemma \ref{lem: vanishing}, where $m$ is some divisible positive integer. 
This implies that $\pi_*\mathcal{O}_X(m(K_X+\Delta)) \rightarrow \pi_*(\mathcal{O}_X(m(K_X+\Delta))\otimes\mathcal{O}_X/\mathcal{J}_1)
$ is surjective.
We consider the commutative diagram
$$
\xymatrix{
\pi^*\pi_*(\mathcal{O}_X(m(K_X+\Delta))\otimes\mathcal{O}_X/\mathcal{J}_1) \ar[r]&
\mathcal{O}_X(m(K_X+\Delta))\otimes\mathcal{O}_X/\mathcal{J}_1
\\
\pi^*\pi_*\mathcal{O}_X(m(K_X+\Delta))\ar[r] \ar[u]&
\mathcal{O}_X(m(K_X+\Delta))\ar[u]
}
$$
We note that $\pi^*\pi_*\mathcal{O}_X(m(K_X+\Delta)) \rightarrow \pi^*\pi_*(\mathcal{O}_X(m(K_X+\Delta))\otimes\mathcal{O}_X/\mathcal{J}_1)
$ is surjective because $\pi^*$ is a right exact functor.
By the definition of $\mathcal{J}_1$, the point $Q$ is isolated in $\operatorname{Supp}(\mathcal{O}_X/\mathcal{J}_1)$. Then $\pi^*\pi_*(\mathcal{O}_X(m(K_X+\Delta))\otimes\mathcal{O}_X/\mathcal{J}_1) \rightarrow
(\mathcal{O}_X(m(K_X+\Delta))\otimes\mathcal{O}_X/\mathcal{J}_1)$ is surjective at $Q$. Therefore, $Q\notin T$.
\end{proof}

\step \label{step: T cap A is empty}
$T\cap A= \emptyset$.

\begin{proof}[Proof of Step \ref{step: T cap A is empty}]
Let $A=\sum_jA_j$ be the connected decomposition.
Now, $\pi(A_j)=P$ for each $j$.
Let $C_i$ be an irreducible component of $A_j$. Then $[C_i]\in N_1(X/Y; W)$ generates an extremal ray of $\overline{\textnormal{NE}}(X/Y; W)$ because $C_i^2<0$. 
By assumption, $C_i$ and $\Delta$ are $\mathbb{Q}$-Cartier.
If $A_j$ is reducible, or it is not isolated in $\operatorname{Supp}\Delta$, then $C_i\cdot (\Delta-C_i)>0$. So, 
$$
C_i\cdot(K_X+C_i)<C_i\cdot(K_X+\Delta)=0.
$$
We note that $\rho(X/Y; W)<\infty$ (cf.~Nakayama's finiteness, \cite[Theorem 4.7]{Fjn22}). By Lemma \ref{Lemma 3.1}, after shrinking $Y$ around $W$, there exists a contraction morphism $\varphi_{C_i}\colon X\rightarrow X'$ over $Y$. 
Let $\pi'\colon X'\rightarrow Y$ be the induced morphism and $\Delta'={\varphi_{C_i}}_*\Delta$. 
Since $(X, Y, W, \pi, \Delta)$ satisfies the condition ($\bigstar$), so does $(X', Y, W, \pi', \Delta')$ by Lemma \ref{lem: preserve conditions}. Moreover, $K_{X'}$ is $\mathbb{Q}$-Cartier and every prime divisor on $X'$ which is mapped into $W$ is $\mathbb{Q}$-Cartier by Lemma \ref{lem: preserve conditions}.

We will check that we can replace $(X, Y, W, \pi, \Delta)$ by $(X', Y, W, \pi', \Delta')$ in order to prove that $T\cap A= \emptyset$. Since $K_X+\Delta$ is $\varphi_{C_i}$-trivial, we have $K_X+\Delta=\varphi_{C_i}^*(K_{X'}+\Delta')$. Then $K_{X'}+\Delta'$ is nef over $Y$.
Let
$$
T'=\operatorname{Supp}\operatorname{Coker}(\pi'^*\pi'_*\mathcal{O}_{X'}(m(K_{X'}+\Delta'))\rightarrow\mathcal{O}_{X'}(m(K_{X'}+\Delta')))\subset X'.
$$
Then, on $X'-T'$, we have the exact sequence
$$
\pi'^*\pi'_*\mathcal{O}_{X'}(m(K_{X'}+\Delta'))\rightarrow\mathcal{O}_{X'}(m(K_{X'}+\Delta'))\rightarrow0.
$$
Taking pull-back by $\varphi_{C_i}\colon X\rightarrow X'$, on $\varphi_{C_i}^{-1}(X'-T')=X-\varphi_{C_i}^{-1}T'$, we have the exact sequence
$$
\varphi_{C_i}^*\pi'^*\pi'_*\mathcal{O}_{X'}(m(K_{X'}+\Delta'))\rightarrow\varphi_{C_i}^*\mathcal{O}_{X'}(m(K_{X'}+\Delta'))\rightarrow0.
$$
Since ${\varphi_{C_i}}_*\varphi_{C_i}^*\mathcal{O}_{X'}(m(K_{X'}+\Delta'))=\mathcal{O}_{X'}(m(K_{X'}+\Delta'))$, on $X-\varphi_{C_i}^{-1}T'$, we have the exact sequence
$$
\varphi_{C_i}^*\pi'^*\pi'_*{\varphi_{C_i}}_*\varphi_{C_i}^*\mathcal{O}_{X'}(m(K_{X'}+\Delta'))\rightarrow\varphi_{C_i}^*\mathcal{O}_{X'}(m(K_{X'}+\Delta'))\rightarrow0.
$$
Then, on $X-\varphi_{C_i}^{-1}T'$, we have the exact sequence
$$
\pi^*\pi_*\mathcal{O}_X(m(K_X+\Delta))\rightarrow\mathcal{O}_X(m(K_X+\Delta))\rightarrow0.
$$
Then $T\subset\varphi_{C_i}^{-1}T'$. This implies that  we can replace $(X, Y, W, \pi, \Delta)$ by $(X', Y, W, \pi', \Delta')$ in order to prove that $T\cap A= \emptyset$. 

We put
\begin{align*}
A^\sharp&:=\sum_{A_j \cap (\operatorname{Supp}\Delta-A_j)=\emptyset}A_j, \textnormal{ and}\\
A^\flat&:=\sum_{A_j \cap (\operatorname{Supp}\Delta-A_j)\neq\emptyset}A_j.
\end{align*}
Therefore, by replacing $X$ with its contraction repeatedly as far as possible, we may assume that $A_j$ is an irreducible and isolated curve in $\operatorname{Supp}\Delta$ when $A_j\subset A^\sharp$. Moreover, by replacing $X$ with its contraction, $A_j$ is contracted to a point $Q$ in $\operatorname{Supp}\Delta$ when $A_j\subset A^\flat$. This point $Q$ is a zero-dimensional non-klt center of $(X, \Delta)$ such that $Q\notin A$ and it satisfies $Q \notin T$ by Step \ref{step: 0-dim lc centers}.
Then we could assume that $A=A^\sharp$ and each $A_j$ is irreducible from the beginning of the proof of Step \ref{step: T cap A is empty}.

We note that since $A_j$ is contracted by a projective morphism, we can treat $A_j$ as an algebraic curve by Serre’s GAGA principle.
Pick a positive integer $m$ such that $m(K_X+\Delta)$ is Cartier. 
If $A_j$ is $\mathbb{P}^1$, then $\mathcal{O}_{A_j}(m(K_X+\Delta)) \cong \mathcal{O}_{A_j}$ because $A_j \cdot\left(K_X+\Delta\right)=0$ and a line bundle on $\mathbb{P}^1$ is trivial if and only if it has degree zero.
If $A_j \ncong \mathbb{P}^1$, then we obtain $h^1(A_j, \mathcal{O}_{A_j})=p_a(A_j) \neq 0$. Since $A_j$ is projective, it has the dualizing sheaf $\omega_{A_j}$. By Serre duality, we obtain $H^0(A_j, \omega_{A_j}) \neq 0$.  We note that
$$
0 \rightarrow \mathcal{T} \rightarrow \mathcal{O}_X\left(K_X+A_j\right) \otimes \mathcal{O}_{A_j} \rightarrow \omega_{A_j} \rightarrow 0
$$
is exact, where the support of $\mathcal{T}$ is of dimension zero by Lemma \ref{lem: adjunction}. Now, we have the exact sequence
$$
H^0(A_j, \mathcal{O}_X(K_X+A_j) \otimes \mathcal{O}_{A_j})\rightarrow
H^0(A_j, \omega_{A_j})\rightarrow
H^1(A_j, \mathcal{T})=0.
$$
Since $H^0(A_j, \omega_{A_j})\neq0$, $H^0(A_j, \mathcal{O}_X(K_X+A_j) \otimes \mathcal{O}_{A_j})\neq0$. Furthermore, as $A_j$ is isolated, 
$$
H^0(A_j, \mathcal{O}_X(K_X+\Delta) \otimes \mathcal{O}_{A_j})=H^0(A_j, \mathcal{O}_X(K_X+A_j) \otimes \mathcal{O}_{A_j}).
$$
Then we see that $\mathcal{O}_{A_j}(m(K_X+\Delta))\cong\mathcal{O}_{A_j}$ because $A_j \cdot\left(K_X+\Delta\right)=0$. Therefore, we obtain $\mathcal{O}_{A}(m(K_X+\Delta)) \cong \mathcal{O}_{A}$. 

Let $f\colon \widetilde{X}\rightarrow X$ be a resolution such that $K_{\widetilde{X}}+\Delta_{\widetilde{X}}=f^*(K_X+\Delta)$. We may assume that
\begin{enumerate}
    \item[(1)] $f^{-1}(A)$ has simple normal crossing support, and
    \item[(2)]  $\operatorname{Supp} f_*^{-1} \Delta \cup \operatorname{Exc}(f)$ is a simple normal crossing divisor on $Y$.
\end{enumerate}
Let $W_3$ be the union of the irreducible components of $\Delta_{\widetilde{X}}^{=1}$ which are mapped into $A$ by $f$. We write $\Delta_{\widetilde{X}}^{=1}=W_3+W_4$. Then
$$
-W_3-\llcorner\Delta_{\widetilde{X}}^{>1}\lrcorner+\ulcorner-(\Delta_{\widetilde{X}}^{<1})\urcorner-(K_{\widetilde{X}}+\{\Delta_{\widetilde{X}}\}+W_4) \sim_\mathbb{Q}-f^*(K_X+\Delta).
$$
We put
$$
\mathcal{J}_2=f_*\mathcal{O}_{\widetilde{X}}(-W_3-\llcorner\Delta_{\widetilde{X}}^{>1}\lrcorner+\ulcorner-(\Delta_{\widetilde{X}}^{<1})\urcorner) \subset \mathcal{O}_X.
$$
Then we can easily check that
$$
0 \rightarrow \mathcal{J}_2 \rightarrow \mathcal{O}_X(-A) \rightarrow \delta \rightarrow 0
$$
is exact, where $\delta$ is a skyscraper sheaf.
As in the proof of Step \ref{step: 0-dim lc centers}, we have
$$
\mathcal{R}^p\pi_*(\mathcal{O}_X(m(K_X+\Delta)) \otimes \mathcal{J}_2)=0
$$
for every $p>0$, where $m$ is some divisible positive integer.
By the above exact sequence, we obtain
$$
\mathcal{R}^p\pi_*(\mathcal{O}_X(m(K_X+\Delta)) \otimes \mathcal{O}_X(-A))=0
$$
for every $p>0$. Then $\pi_*\mathcal{O}_X(m(K_X+\Delta)) \rightarrow \pi_*\mathcal{O}_A(m(K_X+\Delta))
$ is surjective.
We consider the commutative diagram
$$
\xymatrix{
\pi^*\pi_*\mathcal{O}_A(m(K_X+\Delta)) \ar[r]&
\mathcal{O}_A(m(K_X+\Delta))
\\
\pi^*\pi_*\mathcal{O}_X(m(K_X+\Delta))\ar[r] \ar[u]&
\mathcal{O}_X(m(K_X+\Delta))\ar[u]
}
$$
We note that $\pi^*\pi_*\mathcal{O}_X(m(K_X+\Delta)) \rightarrow \pi^*\pi_*\mathcal{O}_A(m(K_X+\Delta))
$ is surjective because $\pi^*$ is a right exact functor, and $\pi^*\pi_*\mathcal{O}_A(m(K_X+\Delta)) \rightarrow
\mathcal{O}_A(m(K_X+\Delta))$ is surjective because $\mathcal{O}_A(m(K_X+\Delta)) \simeq \mathcal{O}_A$.
Therefore, we have $T\cap A=\emptyset$
\end{proof}

\step \label{step: irr comp of B}
Let $E_i$ be an irreducible component of $B$. Then $E_i \not \subset T$.

\begin{proof}[Proof of Step \ref{step: irr comp of B}]
We may assume that $E_i \cap A=\emptyset$ by Step \ref{step: T cap A is empty}. 
When $\pi(E_i)$ is a point, $E_i\cdot(K_X+\Delta)>0$ and so $\mathcal{O}_{E_i}(m(K_X+\Delta))$ is ample. Then $\pi^*\pi_*\mathcal{O}_{E_i}(m(K_X+\Delta))\rightarrow\mathcal{O}_{E_i}(m(K_X+\Delta))$ is surjective.
When $\pi(E_i)$ is a curve, $\pi|_{E_i}$ is a finite morphism. Then $\pi^*\pi_*\mathcal{O}_{E_i}(m(K_X+\Delta))\rightarrow\mathcal{O}_{E_i}(m(K_X+\Delta))$ is surjective.

Let $f\colon \widetilde{X}\rightarrow X$ be the resolution as in the proof of Step \ref{step: T cap A is empty}. We can further assume that
\begin{enumerate}
    \item[(4)] $f^{-1}\left(E_i\right)$ has simple normal crossing support.
\end{enumerate}
Let $W_5$ be the union of the irreducible components of $\Delta_{\widetilde{X}}^{=1}$ which are mapped into $A \amalg E_i$ by $f$. We put $\Delta_{\widetilde{X}}^{=1}=W_5+W_6$. Then
$$
-W_5-\llcorner\Delta_{\widetilde{X}}^{>1}\lrcorner+\ulcorner-(\Delta_{\widetilde{X}}^{<1})\urcorner-(K_{\widetilde{X}}+\{\Delta_{\widetilde{X}}\}+W_6) \sim_\mathbb{Q}-f^*(K_X+\Delta).
$$
We put
$$
\mathcal{J}_3=f_* \mathcal{O}_{\widetilde{X}}(-W_5-\llcorner\Delta_{\widetilde{X}}^{>1}\lrcorner+\ulcorner-(\Delta_{\widetilde{X}}^{<1})\urcorner) \subset \mathcal{O}_X .
$$
Then, as in the proof of Step \ref{step: 0-dim lc centers}, we have
$$
R^p\pi_*(\mathcal{O}_X(m(K_X+\Delta)) \otimes \mathcal{J}_3)=0
$$
for every $p>0$ by Lemma \ref{lem: vanishing}, where $m$ is some divisible positive integer. We note that there exists a short exact sequence
$$
0 \rightarrow \mathcal{J}_3 \rightarrow \mathcal{O}_X\left(-A-E_i\right) \rightarrow \delta^{\prime} \rightarrow 0,
$$
where $\delta^{\prime}$ is a skyscraper sheaf on $X$. Thus,
$$
R^p\pi_*(\mathcal{O}_X(m(K_X+\Delta)) \otimes \mathcal{O}_X(-A-E_i))=0
$$
for every $p>0$. This implies that
$$
\pi_*\mathcal{O}_X(m(K_X+\Delta)) \rightarrow \pi_*\mathcal{O}_{E_i}(m(K_X+\Delta))
$$
is surjective since $\operatorname{Supp} E_i \cap \operatorname{Supp} A=\emptyset$.
Therefore, $E_i \not \subset T$ for every irreducible component $E_i$ of $B$.
\end{proof}

We have checked that $T$ contains no non-klt centers of $(X, \Delta)$, after shrinking $Y$ around $W$ suitably. In fact, it contains no non-klt centers which are components of $\llcorner\Delta\lrcorner$ by Steps \ref{step: T cap A is empty} and \ref{step: irr comp of B}, and it contains no non-klt centers of dimension zero by Steps \ref{step: 0-dim lc centers} and \ref{step: T cap A is empty}. 

Finally, we will prove that $K_X+\Delta$ is semi-ample over a neighborhood of $W$. 

\step
If $T\cap\pi^{-1}(P)=\emptyset$, then there is nothing to prove because $T=\emptyset$ over a neighborhood of $W$. So, we assume that $T\cap\pi^{-1}(P)\neq\emptyset$. 

Since $W=\{P\}$, we only have to prove the result in the case where $Y$ is a Stein space.
In this case, we have
\begin{align*}
T&:=\operatorname{Supp}\operatorname{Coker}(\pi^*\pi_*\mathcal{O}_X(m(K_X+\Delta))\rightarrow\mathcal{O}_X(m(K_X+\Delta)))\\
&=\bigcap_{m(K_X+\Delta)\sim E\geq0}\operatorname{Supp}E 
\end{align*}
set-theoretically by \cite[Lemma 2.7]{EH24}.
By the argument from Step \ref{step: 0-dim lc centers} to Step \ref{step: irr comp of B}, after shrinking $Y$ around $W$ suitably, $T$ is an analytic subset of $X$ and contains no non-klt centers of $(X,\Delta)$. 
We take general members $\Xi_1, \Xi_2, \Xi_3 \in\{E\geq0|m(K_X+\Delta)\sim E\}$ and put $\Theta=\Xi_1+\Xi_2+\Xi_3$. Then $\Theta$ contains no non-klt centers of $(X, \Delta)$ and $K_X+\Delta+\Theta$ is not log canonical at $\operatorname{Supp}\Xi_1\cap\operatorname{Supp}\Xi_2\cap\operatorname{Supp}\Xi_3$ by \cite[Lemma 2.1.3]{Fjn24}.
Since $T\cap\pi^{-1}(P)\neq\emptyset$, $K_X+\Delta+\Theta$ is not log canonical at $\pi^{-1}(P)$.
We put
$$
c=\max \{t\in\mathbb{R}|K_X+\Delta+t\Theta \text { is log canonical at } \pi^{-1}(P)\backslash\operatorname{Nlc}(X, \Delta)\}.
$$
Then we can easily check that $c\in\mathbb{Q}$ and $0<c<1$. In this case,
$$
K_X+\Delta+c\Theta \sim_{\mathbb{Q}}(1+c m)\left(K_X+\Delta\right)
$$
and there exists an lc center $\mathcal{C}$ of $(X, \Delta+c \Theta)$ which is contained in $\cap_{m(K_X+\Delta)\sim E\geq0}\operatorname{Supp}E$, and satisfies $\mathcal{C}\cap \pi^{-1}(P)\neq\emptyset$. We take positive integers $l$ and $n$ such that
$$
l(K_X+\Delta+c\Theta) \sim n m(K_X+\Delta)
$$
Replace $m(K_X+\Delta)$ with $l(K_X+\Delta+c\Theta)$ and apply the previous arguments. Then, after shrinking $Y$ around $W$, $\cap_{kl(K_X+\Delta+c\Theta)\sim E\geq0}\operatorname{Supp}E$ contains no non-klt centers of $(X, \Delta+c\Theta)$ for some positive integer $k$. Since $\mathcal{C}\cap \pi^{-1}(P)\neq\emptyset$, $\mathcal{C}$ is still a non-empty subset of $X$. We obtain $\mathcal{C}\not \subset \cap_{kl(K_X+\Delta+c\Theta)\sim E\geq0}\operatorname{Supp}E$. Therefore,
$$
\bigcap_{knm(K_X+\Delta)\sim E\geq0}\operatorname{Supp}E \subsetneq \bigcap_{m(K_X+\Delta)\sim E\geq0}\operatorname{Supp}E.
$$
This is because there is an lc center $\mathcal{C}$ of $(X, \Delta+c \Theta)$ such that $\mathcal{C} \subset \cap_{m(K_X+\Delta)\sim E\geq0}\operatorname{Supp}E$, and $l(K_X+\Delta+c\Theta)\sim nm(K_X+\Delta)$. By noetherian induction, we deduce that $K_X+\Delta$ is semi-ample over a neighborhood of $W$.

\medskip
The proof of Theorem \ref{Theorem 4.2} is finished.
\end{proof}

\begin{corollary}[{Abundance for $\mathbb{Q}$-divisors, see \cite[Theorem 7.2]{Fjn12}}] \label{Corollary 4.4}
Assume that $(X, Y, W, \pi, \Delta)$ satisfies the condition \textnormal{($\bigstar$)}, $\Delta$ is a $\mathbb{Q}$-divisor,
and that one of \textnormal{Cases (A)-(D)} holds.
If $K_X+\Delta$ is nef over $W$, then it is semi-ample over a neighborhood of $W$. 
\end{corollary}
\begin{proof}
In Case (C), we replace $X$ with its minimal resolution.
In Case (D), any divisor on an open subset of $X$ is $\mathbb{Q}$-Cartier  (cf.~\cite[Theorem 7.3.2]{Ish18}).
So, in these two cases, we may always assume that X is $\mathbb{Q}$-factorial over $W$. So, we only consider Case (A) and Case (B). In Case (B), every prime divisor on $X$ which is mapped into $W$ is $\mathbb{Q}$-Cartier and, shrinking $Y$ around $W$, $\Delta$ is $\mathbb{Q}$-Cartier. By assumption of ($\bigstar$), $K_X$ is also $\mathbb{Q}$-Cartier. So, we only consider Case (A).
In this case, for any point $P\in W$, every prime divisor on $X$ which is included in $\pi^{-1}(P)$ is $\mathbb{Q}$-Cartier. 

As in the proof of Theorem \ref{Theorem 4.2}, we may assume that $K_X+\Delta$ is nef over $Y$.
When $\dim Y=0$, the statement is obvious by Serre’s GAGA principle and \cite[Theorem 7.2]{Fjn12}. Then we may also assume that $\dim Y>0$.
If $\dim Y$=2, then $K_X+\Delta$ is always $\pi$-big. In order to prove that $K_X+\Delta$ is semi-ample over a neighborhood of $W$, we can replace $W$ with $P$ and $Y$ with a Stein open neighborhood of $P$. Then by Theorem \ref{Theorem 4.2}, $K_X+\Delta$ is semi-ample over a neighborhood of $W$. 
If $\operatorname{dim}Y$=1 and $K_X+\Delta$ is $\pi$-big, by Theorem \ref{Theorem 4.2}, $K_X+\Delta$ is semi-ample over a neighborhood of $W$. 
We may assume that $\dim Y$=1 and $K_X+\Delta$ is not $\pi$-big.
Taking the Stein factorization of $\pi\colon X\rightarrow Y$, we may assume that $Y$ is smooth and that $\pi$ has connected fibers.

Let $y \in Y$ be any point, $X_y$ the fiber of $\pi\colon X\rightarrow Y$ over $y$. 
Sard's theorem implies that $X_y$ is smooth for an analytically sufficiently general fiber $X_y$ (cf.~\cite[Corollary 1.8]{Uen75}). Since $X_y$ is connected for each point $y \in Y$, $X_y$ is a smooth irreducible curve for an analytically sufficiently general fiber $X_y$.
Replacing $Y$ with a relatively compact open neighborhood of $W$,
we may assume that $\Delta$ has only finitely many irreducible components. Then $X_y$ and $\Delta$ have no common components for an analytically sufficiently general fiber $X_y$.

Then $(K_X+\Delta)|_{X_y}=K_{X_y}+\Delta|_{X_y}$. Let $\Delta_{X_y}=\Delta|_{X_y}$. Then $\Delta_{X_y}\geq0$ by assumption.
Since $X_y$ is a smooth irreducible curve, 
\begin{align*}
&h^0(X_y, \mathcal{O}_X(m(K_X+\Delta))|_{X_y})-h^1(X_y, \mathcal{O}_X(m(K_X+\Delta))|_{X_y})\\
&=h^0(X_y, \mathcal{O}_{X_y})-h^1(X_y, \mathcal{O}_{X_y})\\
 &=1-h^1(X_y, \mathcal{O}_{X_y}),
\end{align*}
by $\operatorname{deg}(K_X+\Delta)|_{X_y}=0$.
When $h^1(X_y, \mathcal{O}_{X_y})=0$, we have 
$$
H^0(X_y, \mathcal{O}_X(m(K_X+\Delta))|_{X_y})\neq0.
$$
We may assume that $h^1(X_y, \mathcal{O}_{X_y})\geq1$.
By Serre duality, 
$$
1\leq h^1(X_y, \mathcal{O}_{X_y})=h^0(X_y, \mathcal{O}_{X_y}(K_{X_y})).
$$
Since $\Delta_{X_y}\geq0$, there is an effective $\mathbb{Q}$-divisor $D\neq0$ which is $\mathbb{Q}$-linearly equivalent to $K_{X_y}+\Delta_{X_y}$. Then, for a divisible positive integer $m$, we have 
$$
1\leq h^0(X_y, \mathcal{O}_{X_y}(m(K_{X_y}+\Delta_{X_y})))=h^0(X_y, \mathcal{O}_X(m(K_X+\Delta))|_{X_y}).
$$
This implies that $\pi_*\mathcal{O}_X(m(K_X+\Delta))\neq0$. 
Since $Y$ is Stein, we have 
$$
H^0(Y, \pi_*\mathcal{O}_X(m(K_X+\Delta)))=H^0(X, \mathcal{O}_X(m(K_X+\Delta)))\neq0.
$$
Taking the minimal resolution of $X$, we may assume that $X$ is smooth. We note that $\Delta$ is effective, but it is not necessarily boundary, after replacing $X$ by its minimal resolution.

Since $H^0(X, \mathcal{O}_X(m(K_X+\Delta)))\neq0$, there is an effective Cartier divisor $D$ on $X$ such that
$$
\mathcal{O}_X(m(K_X+\Delta))\cong\mathcal{O}_X(D).
$$
It suffices to prove that $D$ is semi-ample over $Y$.
Since we have already replaced $Y$ with a relatively compact open neighborhood of $W$, $D$ has only finitely many irreducible components.
We can write $D=D_{\textnormal{horiz}}+D_{\textnormal{vert}}$, where $D_{\textnormal{horiz}}$ is the horizontal part of $D$ and $D_{\textnormal{vert}}$ is the vertical part of $D$.

Assume that $D_{\textnormal{horiz}}>0$. Then $D_{\textnormal{horiz}}|_{X_y}$ is ample for any irreducible fiber $X_y$ of $\pi\colon X\rightarrow Y$. This implies that $(K_X+\Delta)|_{X_y}$ is ample for general fibers $X_y$. In this case, $K_X+\Delta$ is $\pi$-big. This is a contradiction.

So, we may assume that $D_{\textnormal{horiz}}=0$. We have
$$
0\leq(K_X+\Delta)\cdot D_{\textnormal{vert}}=D_{\textnormal{vert}}^2\leq0,
$$
by Lemma \ref{lem: Zariski}.
Since $D_{\textnormal{vert}}^2=0$, by Lemma \ref{lem: Zariski}, $D_{\textnormal{vert}}$ is proportional to a pull-back of an effective $\mathbb{Q}$-divisor on $Y$. 
In this case, $D_{\textnormal{vert}}=D$ is $\pi$-semi-ample.
\end{proof}

\section{\texorpdfstring{Abundance for $\mathbb{R}$-divisors}{Abundance for R-divisors}} \label{section 5}

\begin{theorem}[{Abundance for $\mathbb{R}$-divisors, see \cite[Theorem 8.1]{Fjn12}}]\label{thm: abundance for R-divisors}
Assume that $(X, Y, W, \pi, \Delta)$ satisfies the condition \textnormal{($\bigstar$)}
and that one of \textnormal{Cases (A)-(D)} holds.
If $K_X+\Delta$ is nef over $W$, then it is semi-ample over a neighborhood of $W$. 
\end{theorem}
The following proof closely follows the argument of \cite[Theorem 8.1]{Fjn12}.

\begin{proof}
As in the proof of Corollary \ref{Corollary 4.4}, it suffices to consider only Case (A).

\setcounter{stepcounter}{0}
\step
Taking a relatively compact open neighborhood of $W$, we may assume that $\operatorname{Supp}\Delta$ has only finitely many irreducible components.
We put $F=\operatorname{Supp}\Delta$ and consider the real vector space $V=\bigoplus_k\mathbb{R}F_k$, where $F=\sum_kF_k$ is the irreducible decomposition. We put $\mathcal{P}=\{D\in V|D\text{ is boundary and }K_X+D\text{ is }\mathbb{R}\text{-Cartier at }\pi^{-1}W\}$. 
Since $\pi^{-1}W$ is a compact subset, $\mathcal{P}$ is a rational polytope. Let $\{R_\lambda\}_{\lambda\in\Lambda}$ be the set of all extremal rays of $\overline{\text{NE}}(X/Y; W)$ spanned by curves. We put $\mathcal{N}=\{D\in\mathcal{P}|(K_X+D)\cdot R_\lambda\geq0\text{ for every }\lambda\in\Lambda\}.$
We note that every $(K_X+\Delta)$-negative extremal ray is spanned by a rational curve $C$ with $0<-(K_X+\Delta)\cdot C\leq4$ (cf.~\cite[Theorem 1.1.6]{Fjn24}).

\medskip
From Step \ref{step: start_rat'l poly} to Step \ref{step: end_rat'l poly}, we will check that $\mathcal{N}$ is a rational polytope in $\mathcal{P}$. The proof of Steps \ref{step: start_rat'l poly} to \ref{step: end_rat'l poly} is essentially the same as the proof of \cite[Proposition 3.2]{Bir11}. 

\step \label{step: start_rat'l poly}
There is a real number $\alpha>0$ such that for any extremal curve $C$ of $\overline{\text{NE}}(X/Y; W)$, if $(K_X+\Delta)\cdot C>0$, then $(K_X+\Delta)\cdot C>\alpha$.

\begin{proof}[Proof of Step \ref{step: start_rat'l poly}]
If $\Delta$ is a $\mathbb{Q}$-divisor, then the statement is trivially true. If $\Delta$ is not a $\mathbb{Q}$-divisor, let $\Delta_1$, $\cdots$, $\Delta_r$ be the vertices of $\mathcal{P}$ and let $a_1$, $\cdots$, $a_r$ be nonnegative real numbers such that $\Delta=\sum_ia_i\Delta_i$ and $\sum_ia_i=1$. Then, $(K_X+\Delta)\cdot C=\sum_ia_i(K_X+\Delta_i)\cdot C$ and if $(K_X+\Delta)\cdot C<1$, then there are only finitely many possibilities for the intersection numbers $(K_X+\Delta_i)\cdot C$ because $(K_X+\Delta_i)\cdot C\geq-4$. So, the existence of $\alpha$ is clear.   
\end{proof}

\step \label{step: second_rat'l poly}
There is a real number $\delta>0$ such that if $D\in\mathcal{P}$, $\lVert D-\Delta\rVert<\delta$ and $(K_X+D)\cdot R\leq0$ for an extremal ray $R$ of $\overline{\text{NE}}(X/Y; W)$, then $(K_X+\Delta)\cdot R\leq0$. 

\begin{proof}[Proof of Step \ref{step: second_rat'l poly}]
If the statement is not true, then there is an infinite sequence of $D_\lambda\in\mathcal{P}$ and extremal rays $R_\lambda$ of $\overline{\text{NE}}(X/Y; W)$ such that for each $\lambda$ we have $(K_X+D_\lambda)\cdot R_\lambda\leq0$, $(K_X+\Delta)\cdot R_\lambda>0$, and $\lVert D_\lambda-\Delta\rVert$ converges to 0. For each $\lambda$, there are nonnegative real numbers $a_{1, \lambda}$, $\cdots$, $a_{r, \lambda}$ such that $D_\lambda=\sum_ia_{i, \lambda}\Delta_i$ and $\sum_ia_{i, \lambda}=1$. 
The representation of $\Delta=\sum_ia_i\Delta_i$ and $D_\lambda=\sum_ia_{i, \lambda}\Delta_i$ is not necessarily unique now. 
Taking an infinite subsequence of $\{D_\lambda\}$, we may assume that there is some simplex $\mathcal{S}$ which is spanned by some vertices of $\mathcal{P}$ such that every $D_\lambda$ is in the inner of $\mathcal{S}$. 
We can write $D_\lambda=\sum_ia_{i, \lambda}\Delta_i$, where $a_{i, \lambda}>0$ if and only if $\Delta_i$ is a vertex of $\mathcal{S}$, and this representation is unique. Since $\lVert D_\lambda-\Delta\rVert$ converges to 0, the limits $\lim_{\lambda\to\infty}a_{i, \lambda}$ exist and we have $\Delta=\sum_i(\lim_{\lambda\to\infty}a_{i, \lambda})\Delta_i$. We assume that $a_i=\lim_{\lambda\to\infty}a_{i, \lambda}$.
After replacing the sequence with an infinite subsequence, we can assume that the sign of $(K_X+\Delta_i)\cdot R_\lambda$ is independent of $\lambda$, and that for each $\lambda$ we have an extremal curve $C_\lambda$ for $R_\lambda$.
Now, if $(K_X+\Delta_i)\cdot C_\lambda\leq0$, then it is bounded from below, hence there are only finitely many possibilities for this number and we could assume that it is independent of $\lambda$.
On the other hand, if $a_i\neq0$, then $(K_X+\Delta_i)\cdot C_\lambda$ is bounded from below and above because $(K_X+D_\lambda)\cdot C_\lambda=\sum_ia_{i, \lambda}(K_X+\Delta_i)\cdot C_\lambda\leq0$, hence there are only finitely many possibilities for $(K_X+\Delta_i)\cdot C_\lambda$ and we could assume that it is independent of $\lambda$.
Assume that $a_i\neq0$ for $1\leq i\leq l$ but $a_i=0$ for $i>l$. Then, it is clear that
\begin{align*}
(K_X+D_\lambda)\cdot C_\lambda=(K_X+\Delta)\cdot C_\lambda&+\sum_{i\leq l}(a_{i, \lambda}-a_i)(K_X+\Delta_i)\cdot C_\lambda\\
&+\sum_{i>l}a_{i, \lambda}(K_X+\Delta_i)\cdot C_\lambda  
\end{align*}
would be positive by Step \ref{step: start_rat'l poly} if $\lambda\gg0$, which gives a contradiction.
\end{proof}

\step \label{step: end_rat'l poly}
The set $\mathcal{N}$ is a rational polytope.    

\begin{proof}[Proof of Step \ref{step: end_rat'l poly}]
We may assume that for each $\lambda\in\Lambda$ there is some $D\in\mathcal{P}$ such that $(K_X+D)\cdot R_\lambda<0$, in particular, $(K_X+\Delta_i)\cdot R_\lambda<0$ for a vertex $\Delta_i$ of $\mathcal{P}$. Since the set of such extremal rays is discrete, we may assume that $\Lambda\subset\mathbb{N}$. Obviously, $\mathcal{N}$ is a convex compact subset of $\mathcal{P}$. If $\Lambda$ is finite, the claim is trivial. So we may assume that $\Lambda=\mathbb{N}$. By Step \ref{step: second_rat'l poly} and by the compactness of $N$, there are $D_1,\cdots,D_n\in\mathcal{N}$ and $\delta_1,\cdots,\delta_n>0$ such that $\mathcal{N}$ is covered by $\mathscr{B}_i=\{D\in\mathcal{P}|\lVert D-D_i\rVert<\delta_i\}$, and such that if $D\in\mathscr{B}_i$ with $(K_X+D)\cdot R_\lambda<0$ for some $\lambda$, then $(K_X+D_i)\cdot R_\lambda=0$. If $\Lambda_i=\{\lambda\in\Lambda|(K_X+D)\cdot R_\lambda<0\text{ for some }D\in\mathscr{B}_i\}$, then by construction $(K_X+D_i)\cdot R_\lambda=0$ for any $\lambda\in\Lambda_i$. Then, since the $\mathscr{B}_i$ give an open cover of $\mathcal{N}$, we have $\mathcal{N}=\bigcap_{1\leq i\leq l}\mathcal{N}_{\Lambda_i}$, where $\mathcal{N}_{\Lambda_i}=\{D\in\mathcal{P}|(K_X+D)\cdot R_\lambda\geq0\text{ for all }\lambda\in\Lambda_i\}$. So, it is enough to prove that each $\mathcal{N}_{\Lambda_i}$ is a rational polytope and replacing $\Lambda$ with $\Lambda_i$, we could assume from the beginning that there is some $D\in\mathcal{N}$ such that $(K_X+D)\cdot R_\lambda=0$ for every $\lambda\in\Lambda$. If $\operatorname{dim}\mathcal{P}=1$, this already proves the claim. If $\operatorname{dim}\mathcal{P}>1$, let $\mathcal{P}^1,\cdots,\mathcal{P}^p$ be the proper faces of $\mathcal{P}$. Then, each $\mathcal{N}^j=\mathcal{N}\cap\mathcal{P}^j$ is a rational polytope by induction. Moreover, for each $D''\in\mathcal{N}$ which is not $D$, there is $D'$ on some proper face of $\mathcal{P}$ such that $D''$ is on the line segment determined by $D$ and $D'$. Since $(K_X+D)\cdot R_\lambda=0$ for every $\lambda\in\Lambda$, if $D'\in\mathcal{P}^j$, then $D'\in\mathcal{N}^j$. Hence $\mathcal{N}$ is the convex hull of $D$ and all the $\mathcal{N}^j$. Now, there is a finite subset $\Lambda'\subset\Lambda$ such that $\cup_j\mathcal{N}^j=\mathcal{N}_{\Lambda'}\cap(\cup_j\mathcal{P}^j)$. But then the convex hull of $D$ and $\mathcal{N}^j$ is just $\mathcal{N}_{\Lambda'}$ and we are done.
\end{proof}
We have checked that $\mathcal{N}$ is a rational polytope in $\mathcal{P}$. 

Finally, we will prove that $K_X+\Delta$ is semi-ample over a neighborhood of $W$. 

\step \label{step: final R-divisor}
We note that $\mathcal{N}=\{D\in\mathcal{P}|K_X+D\text{ is nef}\}$. By the above construction, we have $\Delta\in\mathcal{N}$. Let $\mathcal{F}$ be the minimal face of $\mathcal{N}$ containing $\Delta$. Then we can take $\mathbb{Q}$-divisors $\Delta_1,\cdots,\Delta_q$ on $X$ and positive real numbers $a_1,\cdots,a_q$ such that $\Delta_i$ is in the relative interior of $\mathcal{F}$ for every $i$, $K_X+\Delta=\sum_ia_i(K_X+\Delta_i)$ and $\sum_ia_i=1$. By Corollary \ref{Corollary 4.4}, $K_X+\Delta_i$ is semi-ample over a neighborhood of $W$ for every $i$ since $K_X+\Delta_i$ is nef over $W$. Therefore, $K_X+\Delta$ is semi-ample over a neighborhood of $W$. 
\end{proof}

\section{Proof of the main theorem and corollaries} \label{section 6}

\begin{proof}[Proof of Theorem \ref{Theorem 1.1}]
Theorems \ref{thm: MMP} and \ref{thm: abundance for R-divisors} imply Theorem \ref{Theorem 1.1}.
\end{proof}

\begin{proof}[Proof of Corollary \ref{Corollary 1.5}]
We may assume that $\pi_*\mathcal{O}_X(\llcorner m(K_X+\Delta)\lrcorner)\neq0$ for some $m\geq0$. In this case, $K_X+\Delta$ is $\pi$-pseudo-effective. By Theorem \ref{Theorem 1.1}, shrinking $Y$ around $W$, we may assume that $K_X+\Delta$ is $\pi$-semi-ample. Then Corollary \ref{Corollary 1.5} follows from \cite[2.36]{Fjn22}.
\end{proof}

\begin{proof}[Proof of Corollary \ref{Corollary 1.6}]
In Case (C),  we replace $X$ with its minimal resolution. So, we may always assume that X has only rational singularities. Then we only have to prove the statement in Case (D).
In this case, any divisor on an open subset of $X$ is $\mathbb{Q}$-Cartier (cf.~\cite[Theorem 7.3.2]{Ish18}). Then, for each point $y\in Y$, $X$ is $\mathbb{Q}$-factorial over $W=\{y\}$. We note that $\rho(X/Y; W)<\infty$ (cf.~Nakayama's finiteness, \cite[Theorem 4.7]{Fjn22}).
By Theorem \ref{Theorem 1.1}, if $K_X+\Delta$ is $\pi$-pseudo-effective, we may assume that $K_X+\Delta$ is $\pi$-semi-ample over some open neighborhood of $y$. By \cite[2.36]{Fjn22}, the log canonical ring $R(X,\Delta)$ is locally finitely generated over some open neighborhood of $y$. Then it is locally finitely generated on $Y$.
\end{proof}

\bibliographystyle{amsalpha} 
\bibliography{ref}

\providecommand{\bysame}{\leavevmode\hbox to3em{\hrulefill}\thinspace}
\providecommand{\MR}{\relax\ifhmode\unskip\space\fi MR }
\providecommand{\MRhref}[2]{%
  \href{http://www.ams.org/mathscinet-getitem?mr=#1}{#2}
}
\providecommand{\href}[2]{#2}
\begin{thebibliography}{BHPV04}

\bibitem[BHPV04]{BHPV04}
Wolf~P. Barth, Klaus Hulek, Chris A.~M. Peters, and Antonius {Van de Ven}, \emph{Compact complex surfaces}, second ed., Ergebnisse der Mathematik und ihrer Grenzgebiete. 3. Folge. A Series of Modern Surveys in Mathematics [Results in Mathematics and Related Areas. 3rd Series. A Series of Modern Surveys in Mathematics], vol.~4, Springer-Verlag, Berlin, 2004. \MR{2030225}

\bibitem[Bir11]{Bir11}
Caucher Birkar, \emph{On existence of log minimal models {II}}, Journal für die reine und angewandte Mathematik \textbf{2011} (2011), no.~658, 99--113.

\bibitem[EH24]{EH24}
Makoto Enokizono and Kenta Hashizume, \emph{Minimal model program for log canonical pairs on complex analytic spaces}, 2024, \href{https://arxiv.org/abs/2404.05126}{arXiv:2404.05126}.

\bibitem[FM25]{FM25}
Osamu Fujino and Nao Moriyama, \emph{On a vanishing theorem for surfaces}, 2025, \href{https://arxiv.org/abs/2511.19783}{arXiv:2511.19783}.

\bibitem[Fuj12]{Fjn12}
Osamu Fujino, \emph{Minimal model theory for log surfaces}, Publ. Res. Inst. Math. Sci. \textbf{48} (2012), no.~2, 339--371. \MR{2928144}

\bibitem[Fuj22]{Fjn22}
\bysame, \emph{Minimal model program for projective morphisms between complex analytic spaces}, 2022, \href{https://arxiv.org/abs/2201.11315}{arXiv:2201.11315}.

\bibitem[Fuj24]{Fjn24}
\bysame, \emph{Cone and contraction theorem for projective morphisms between complex analytic spaces}, MSJ Memoirs, vol.~42, Mathematical Society of Japan, Tokyo, 2024. \MR{4845344}

\bibitem[Fuj25]{Fjn23}
\bysame, \emph{Vanishing theorems for projective morphisms between complex analytic spaces}, Math. Res. Lett. \textbf{32} (2025), no.~3, 739--769. \MR{4945100}

\bibitem[Has16]{Hsh16}
Kenta Hashizume, \emph{Finite generation of adjoint ring for log surfaces}, J. Math. Sci. Univ. Tokyo \textbf{23} (2016), no.~4, 741--761. \MR{3588261}

\bibitem[Ish18]{Ish18}
Shihoko Ishii, \emph{Introduction to singularities}, second ed., Springer, Tokyo, 2018. \MR{3838338}

\bibitem[KM98]{KM98}
J\'anos Koll\'ar and Shigefumi Mori, \emph{Birational geometry of algebraic varieties}, Cambridge Tracts in Mathematics, vol. 134, Cambridge University Press, Cambridge, 1998, With the collaboration of C. H. Clemens and A. Corti, Translated from the 1998 Japanese original. \MR{1658959}

\bibitem[Mat02]{Mts02}
Kenji Matsuki, \emph{Introduction to the {M}ori program}, Universitext, Springer-Verlag, New York, 2002. \MR{1875410}

\bibitem[RRV71]{RRV71}
Jean~P. Ramis, Gabriel Ruget, and Jean~L. Verdier, \emph{Dualit\'e{} relative en g\'eom\'etrie analytique complexe}, Invent. Math. \textbf{13} (1971), 261--283. \MR{308439}

\bibitem[Sak84]{Fum84}
Fumio Sakai, \emph{Weil divisors on normal surfaces}, Duke Math. J. \textbf{51} (1984), no.~4, 877--887. \MR{771385}

\bibitem[Tan14]{Tnk14}
Hiromu Tanaka, \emph{Minimal models and abundance for positive characteristic log surfaces}, Nagoya Math. J. \textbf{216} (2014), 1--70. \MR{3319838}

\bibitem[Uen75]{Uen75}
Kenji Ueno, \emph{Classification theory of algebraic varieties and compact complex spaces}, Lecture Notes in Mathematics, vol. Vol. 439, Springer-Verlag, Berlin-New York, 1975, Notes written in collaboration with P. Cherenack. \MR{506253}

\end{thebibliography}

\end{document}